\documentclass[12pt]{amsart}
\usepackage[T2A]{fontenc}
\usepackage[english]{babel}
\usepackage{amsaddr}
\usepackage{amsmath}
\usepackage{amssymb}
\usepackage{amsfonts}
\usepackage{epsfig}
\usepackage{srcltx}
\usepackage{subfigure}
\usepackage{cite}
\usepackage{float}
\usepackage{mathtools}
\usepackage[a4paper, mag=1000, includefoot, left=2cm, right=2cm, top=2cm, bottom=2cm, headsep=1cm, footskip=1cm]{geometry}
\usepackage[unicode,colorlinks]{hyperref}
\hypersetup{colorlinks=true, citecolor=blue, linkcolor=blue}

\newtheorem{Th}{Theorem}

\newtheorem{Def}{Definition}
\begin{document}
\thispagestyle{empty}

\title[Bifurcations in asymptotically autonomous systems]{Bifurcations in asymptotically autonomous Hamiltonian systems under multiplicative noise}
\author[O.A. Sultanov]{Oskar A. Sultanov}

\address{
Chebyshev Laboratory, St. Petersburg State University, 14th Line V.O., 29, Saint Petersburg 199178 Russia.}
\email{oasultanov@gmail.com}


\maketitle
{\small

{\small
\begin{quote}
\noindent{\bf Abstract.} The effect of multiplicative stochastic perturbations on Hamiltonian systems on the plane is investigated. It is assumed that perturbations fade with time and preserve a stable equilibrium of the limiting system. The paper investigates bifurcations associated with changes in the stability of the equilibrium and with the appearance of new stochastically stable states in the perturbed system. It is shown that depending on the structure and the parameters of the decaying perturbations the equilibrium can remain stable or become unstable. In some intermediate cases, a practical stability of the equilibrium with estimates for the length of the stability interval is justified. The proposed stability analysis is based on a combination of the averaging method and the construction of stochastic Lyapunov functions.

\medskip



\end{quote}
}

\section{Introduction}
The influence of perturbations on dynamical systems is a classical problem with a wide range of applications. For deterministic autonomous systems, the solution of such a problem is effectively solved by the theory of stability and bifurcations~\cite{GH83,GL94}. The effect of perturbations on asymptotically autonomous systems has been studied in many papers. Under some conditions, the solutions of perturbed systems can retain the asymptotic properties of trajectories of the corresponding limiting autonomous systems~\cite{RB53,LM56,LDP74}. In general this is not the case: behavior of perturbed and unperturbed trajectories can differ significantly~\cite{HRT94}. The qualitative and asymptotic properties of solutions to such systems depend both on the properties of the unperturbed system and on the structure of disturbances. Bifurcation phenomena in asymptotically autonomous systems were discussed in~\cite{LRS02,KS05,MR08,OS21DCDS}.

This paper focuses on the stochastic perturbations of dynamical systems. It is well known that even weak random disturbances can lead to significant changes in the behavior of trajectories~\cite{FW98}. See, for example,~\cite{RKh64,BG02,AMR08,BKGT08,DNR11,BHW12,KT13,TW15}, where the influence of autonomous stochastic perturbations on qualitative properties of solutions is discussed. Stochastic bifurcations associated with qualitative changes in the profile of stationary probability densities, in the Lyapunov spectrum function or in the dichotomy spectrum were investigated in~\cite{NSN90,CF98,BRS09,BRR15,CDLR17,DELR18} for systems of stochastic differential equations with time-independent coefficients. The influence of damped stochastic perturbations on the long-term behavior of solutions of scalar autonomous systems was discussed in~\cite{AGR09,ACR11}. To the best of the author's knowledge, the effect of decaying stochastic perturbations on bifurcations in asymptotically autonomous systems has not been thoroughly investigated.

The present paper is devoted to studying the effect of multiplicative stochastic perturbations on asymptotically autonomous Hamiltonian systems on the plane. It is assumed that the intensity of noise fades with time and the limiting system has a neutrally stable equilibrium. Bifurcations associated with changes in the stochastic stability of the equilibrium, as well as with the appearance of new stable states, are discussed. 

The paper is organized as follows. In Section~\ref{sec1}, the formulation of the problem is given and the class of fading perturbations is described. The main results are presented in Section~\ref{sec2}. The proofs are contained in the subsequent sections. In Section~\ref{Sec3}, changes of variables are constructed that simplifies the system in the first asymptotic terms at infinity. This transformation consists of a transition to energy-angle variables associated with the parameters of the general solution of the unperturbed Hamiltonian system, and a specific nearly-identity transformation of the energy variable. The study of the structure of the simplified equations and the nonlinear stability analysis lead to a description of possible bifurcations in the system. The stability analysis based on the construction of stochastic Lyapunov functions is contained in sections~\ref{Sec41},~\ref{Sec42} and \ref{Sec43}. In Section~\ref{SecEx}, the proposed theory is applied to the examples of nonlinear systems with decaying stochastic perturbations. The paper concludes with a brief discussion of the results obtained.

\section{Problem statement}\label{sec1}
Consider the system of It\^{o} stochastic differential equations on the plane:
\begin{equation}\label{FulSys}
	d {\bf z} = {\bf b}({\bf z},t) dt + {\bf B}({\bf z},t)d{\bf w}(t), \quad t>s>0, \quad {\bf z}(s)={\bf z}_0\in\mathbb R^2,
\end{equation}
where
\begin{gather*}
{\bf z}=\begin{pmatrix}x\\ y\end{pmatrix}, \quad 
{\bf b}({\bf z},t)\equiv \begin{pmatrix}\partial_y H(x,y,t) \\ -\partial_x H(x,y,t)+F(x,y,t)\end{pmatrix}, \quad
{\bf B}({\bf z},t)\equiv
	\begin{pmatrix}
		B_{1,1}(x,y,t) & B_{1,2}(x,y,t) \\
		B_{2,1}(x,y,t) & B_{2,2}(x,y,t)
\end{pmatrix},
\end{gather*} 
and ${\bf w}(t)=(w_1(t),w_2(t))^T$ is a two-dimensional Wiener process defined on a probability space $(\Omega,\mathcal F,\mathbb P)$. The functions $H(x,y,t)$, $F(x,y,t)$, $B_{i,j}(x,y,t)$  are defined for all $(x,y,t)\in\mathbb R^2\times \mathbb R_+$, are infinitely differentiable and do not depend on $\omega\in\Omega$. It is assumed that 
\begin{gather}\label{zero}
{\bf b}(0,t)\equiv 0, \quad {\bf B}(0,t)\equiv 0,
\end{gather}
and there exists $M>0$ such that
\begin{gather}\label{lip}
|{\bf b}({\bf z}_1,t)-{\bf b}({\bf z}_2,t)|\leq M |{\bf z}_1-{\bf z}_2|, \quad \|{\bf B}({\bf z}_1,t)-{\bf B}({\bf z}_2,t)\|\leq M |{\bf z}_1-{\bf z}_2|
\end{gather}
for all ${\bf z}_1, {\bf z}_2\in \mathbb R^2$ and $t\geq s$, where $|{\bf z}|=\sqrt{x^2+y^2}$ and $\|\cdot\|$ is the operator norm coordinated with the norm $|\cdot|$ of $\mathbb R^2$. These constraints on the coefficients guarantee the existence and uniqueness of a continuous (with probability one) solution ${\bf z}(t)=(x(t),y(t))^T$ for all $t\geq s$ for any initial point ${\bf z}_0\in\mathbb R^2$ (see, for example, \cite[Sec. 5.2]{BOks98}).

Furthermore, it is assumed that system \eqref{FulSys} is asymptotically autonomous, and for every compact $\mathcal D\subset\mathbb R^2$  
\begin{equation*}
	\lim_{t\to\infty} H(x,y,t)=H_0(x,y), \quad \lim_{t\to\infty} F(x,y,t)=\lim_{t\to\infty} B_{i,j}(x,y,t)=0 
\end{equation*}
for all $(x,y)\in \mathcal D$ and $i,j\in\{1,2\}$. The limiting autonomous system
\begin{equation}\label{LimSys}
	\frac{dx}{dt}=\partial_y H_0(x,y), \quad \frac{dy}{dt}=-\partial_x H_0(x,y)
\end{equation}
is assumed to have the isolated fixed point $(0,0)$ of center type. Without loss of generality, it is assumed that 
\begin{equation}\label{H0as}
	H_0(x,y)=\frac{|{\bf z}|^2}{2}+\mathcal O(|{\bf z}|^3), \quad |{\bf z}|\to 0,
\end{equation}
and there exist $E_0>0$ and $r>0$ such that for all $E\in [0,E_0]$ the level lines $\{(x,y)\in\mathbb R^2: H_0(x,y)=E\}$, lying in $\mathcal B_r=\{(x,y)\in\mathbb R^2: |{\bf z}|\leq r\}$, define a family of closed curves on the phase space $(x,y)$ parametrized by the parameter $E$. 
Each of these curves corresponds to a periodic solution $\hat x(t,E)$, $\hat y(t,E)$ of system \eqref{LimSys} with  
a period $T(E)=2\pi/\nu(E)$, where $\nu(E)\neq 0$ for all $E\in[0,E_\ast]$ and $\nu(E)=1+\mathcal O(E)$ as $E\to 0$.  The value $E=0$ corresponds to the fixed point $(0,0)$. It is also assumed that $\mathcal B_r$ does not contain any fixed points of the limiting system, except for the origin.

Damped perturbations of the limiting system \eqref{LimSys} are described by functions with power-law asymptotic expansions:
\begin{equation}\label{HFBas}
	\begin{split}
		  H(x,y,t)=H_0(x,y)+\sum_{k=1}^\infty t^{-\frac{k}{q}} H_k(x,y), \quad
		 F(x,y,t)=\sum_{k=1}^\infty t^{-\frac{k}{q}} F_k(x,y), \\ 
		  B_{i,j}(x,y,t)=\sum_{k=1}^\infty t^{-\frac{k}{q}} B_{i,j,k}(x,y),  \quad t\to\infty
	\end{split}
\end{equation}
with $q\in\mathbb Z_+$.  Note that the series in \eqref{HFBas} are assumed to be asymptotic as $t\to\infty$ uniformly for all $(x,y)\in \mathcal B_r$ (see, for example, \cite[\S 1]{MVF89}). Such decaying perturbations appear, for example, in the study of Painlev\'{e} equations~\cite{IKNF06,BG08}, resonance and phase-locking phenomena~\cite{LK14,OS21} and in many other problems associated with nonlinear and non-autonomous systems~\cite{KF13,CYZZ18,CH21}.

It can easily be checked that the rational powers of the form $k/q$ with $q>1$ in \eqref{HFBas} can be reduced to the integer exponents $k$ by the change of the time variable $\theta=t^{1/q}$ in system \eqref{FulSys}. In this case, the growing factor $\theta^{q-1}$ appears in the drift term of system \eqref{FulSys}. This indicates that the problem of a long-term behavior of solutions is singularly perturbed~\cite{FV05}: a global behavior of solutions cannot be derived from a corresponding limiting equations. Note also that the decaying factors $t^{-k/q}$ in the right-hand side of system \eqref{FulSys} can be rewritten in the form $\epsilon^{k/q}\varsigma^{-k/q}$ with a new independent variable $\varsigma=\epsilon t$ and a small parameter $0<\epsilon\ll1$. Although in some cases the asymptotic solution of such a problem as $\epsilon\to0$ and $\varsigma=\mathcal O(1)$ can give a long-term approximation for solutions in the original variable $t=\varsigma/\epsilon$, this approach is usually not used when investigating the behavior of solutions at infinity~\cite{WW66,FO74,BAD04,IKNF06,LK09,ACR11}. Moreover, in some cases such method is inapplicable due to the appearance of metastable states under decaying deterministic perturbations (see, for example,~\cite{OSIJBC21}). In the present paper, a small parameter is not introduced.

The simplest example is given by the following linear system:
\begin{equation}\label{ex0}
\begin{split}
	 dx=ydt,\quad
	 dy=(-x+t^{-1}\lambda y)\,dt+ t^{-\frac 12} \mu x\,dw_2(t), \quad t>1,
\end{split}
\end{equation}
with $\lambda,\mu={\hbox{\rm const}}$. This system is of the form \eqref{FulSys} with $q=2$, $H(x,y,t)\equiv H_0(x,y)\equiv |{\bf z}|^2/2$, $F(x,y,t)\equiv t^{-1}\lambda y$, $B_{1,1}(x,y,t)\equiv B_{1,2}(x,y,t)\equiv B_{2,1}(x,y,t)\equiv 0$ and $B_{2,2}(x,y,t)\equiv t^{-1/2} \mu x$. It can easily be checked that the autonomous system with $\lambda=0$ and $\mu=0$ has a periodic general solution $x_\ast(t;E,\varphi)=\sqrt{2E}\cos (t+\varphi)$, $y_\ast(t;E,\varphi)=-\sqrt{2E}\sin (t+\varphi)$ with a period $T(E)\equiv 2\pi$. The asymptotics for a two-parameter family of solutions to the perturbed deterministic system with $\lambda\neq 0$ and $\mu=0$ is constructed with the WKB approximations~\cite{WW66}: 
\begin{gather*}
	x(t)=t^{\frac{\lambda}{2}}\left(x_\ast(t;E,\varphi)+\mathcal O(t^{-1})\right),\quad 
	y(t)=t^{\frac{\lambda}{2}}\left(y_\ast(t;E,\varphi)+\mathcal O(t^{-1})\right),\quad 
	t\to\infty.
\end{gather*}
In this case, the stability of the equilibrium $(0,0)$ depends on the sign of the parameter $\lambda$ (see Fig.~\ref{Fig12}, a). Numerical analysis of system \eqref{ex0} with $\lambda\neq 0$ and $\mu\neq 0$ shows that decaying stochastic perturbations can lead to the shift of the stability boundary: the stability of the equilibrium $(0,0)$ changes as $\lambda$ passes through a certain critical value $\lambda_\ast(\mu)$ (see Fig.~\ref{Fig12}, b). More tricky examples are considered in Section~\ref{SecEx}.
\begin{figure}
\centering
\subfigure[$\mu=0$ ]{\includegraphics[width=0.4\linewidth]{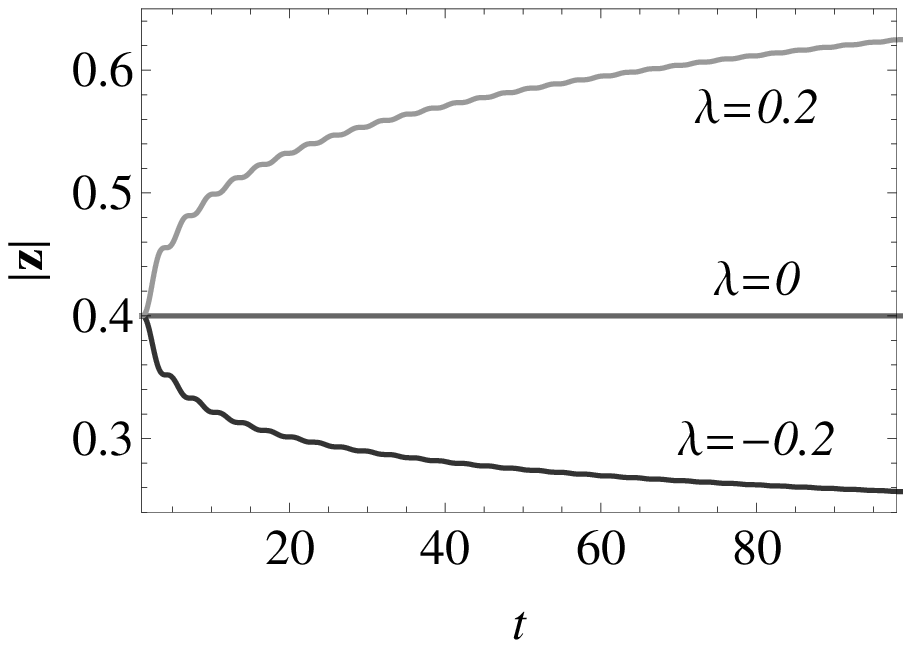}}
\hspace{4ex}
\subfigure[$\mu=1$]{\includegraphics[width=0.4\linewidth]{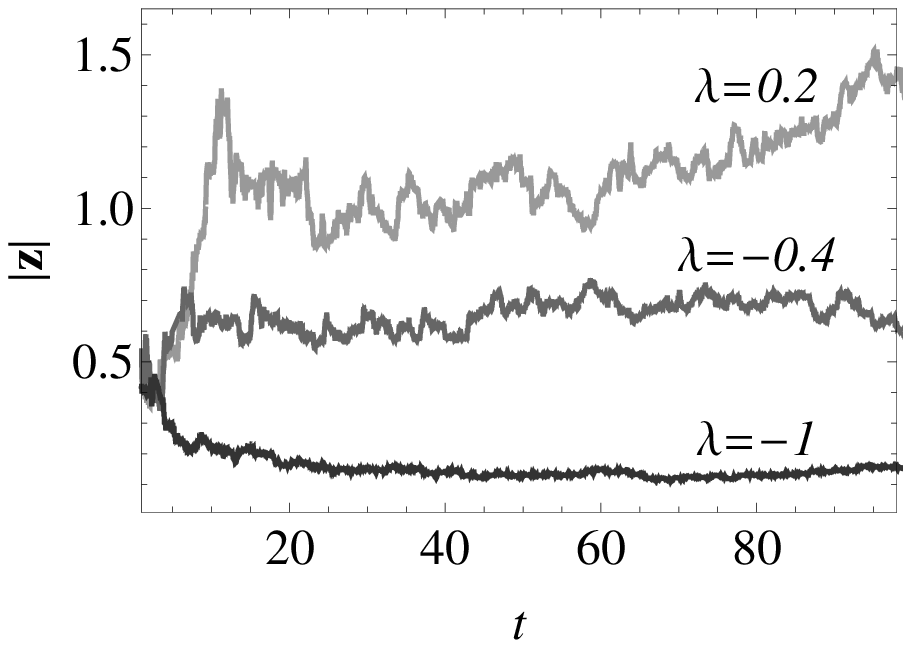}}
\caption{\small The evolution of $|{\bf z}(t)|=\sqrt{x^2(t)+y^2(t)}$ for sample paths of the solutions to system \eqref{ex0} with initial data $x(1)=0.4$, $y(1)=0$.} \label{Fig12}
\end{figure}

In the general case, nonlinear terms of equations with decaying parameters also affect the qualitative behavior of solutions to systems of the form \eqref{FulSys}. The goal of this paper is to describe the stability conditions for system \eqref{FulSys} and to reveal the role of decaying perturbations in the corresponding local bifurcations, associated with changes in stochastic stability of the solution ${\bf z}(t)\equiv 0$.

\section{Main results}\label{sec2}

First, note that system \eqref{FulSys} can be simplified at least in the first terms of the asymptotics, by a suitable transformation of the variables. 

Let $\hat x(t,E)$, $\hat y(t,E)$ be one of $T(E)$-periodic solutions of system \eqref{LimSys} such that 
$H_0(\hat x(t,E),\hat y(t,E))\equiv E$ as $E\in[0,E_0]$. Define $\mathcal D(E_0)=\{(x,y)\in\mathcal B_r: 0\leq H_0(x,y)\leq E_0\}$. Then, we have the following.

\begin{Th}\label{Th1}
Let system \eqref{FulSys} satisfy \eqref{zero}, \eqref{lip}, \eqref{H0as} and \eqref{HFBas}. 
Then for all $N\geq 1$ there exist $t_\ast\geq s$ and the chain of transformations $(x,y)\to (E,\varphi)\to (v,\varphi)$, defined by
\begin{gather}
\label{exch1} x(t)=\hat x\left(\frac{\varphi(t)}{\nu(E(t))},E(t)\right), \quad y(t)=\hat y\left(\frac{\varphi(t)}{\nu(E(t))},E(t)\right), \\
\label{exch11} v(t)=V_N(E(t),\varphi(t),t), \quad V_N(E,\varphi,t)\equiv E+\sum_{k=1}^N t^{-\frac{k}{q}} v_k(E,\varphi),
\end{gather}
such that for all $(x,y)\in \mathcal D(E_0)$ and $t\geq t_\ast$ system \eqref{FulSys} can be transformed into
\begin{gather}\label{Veq}
\begin{split}
dv=&\left( \sum_{k=1}^N t^{-\frac{k}{q}}\Lambda_k(v)+\tilde \Lambda_N(v,\varphi,t)\right) \,dt+\sum_{j=1}^2 \alpha_{1,j}(v,\varphi,t) \, dw_j(t),\\
{d\varphi} =&\left(\nu(v)+\tilde G_N (v,\varphi,t)\right)dt+ \sum_{j=1}^2 \alpha_{2,j}(v,\varphi,t) dw_j(t),
\end{split}
\end{gather}
where $v_k(E,\varphi)$, $\tilde \Lambda_{N}(v,\varphi,t)$, $\tilde G_N(v,\varphi,t)$ and $\alpha_{i,j}(v,\varphi,t)$ are $2\pi$-periodic in $\varphi$, 
$\Lambda_k(0)\equiv \tilde \Lambda_N(0,\varphi,t)\equiv \alpha_{1,j}(0,\varphi,t)\equiv  v_k(0,\varphi)\equiv 0$,
and the following estimates hold:
\begin{gather}\label{RemEst}
\tilde \Lambda_N(v,\varphi,t)=\mathcal O(t^{-\frac{N+1}{q}}), \quad \tilde G_N(v,\varphi,t)=\mathcal O(t^{-\frac{1}{q}}), \quad \alpha_{i,j}(v,\varphi,t)=\mathcal O(t^{-\frac{1}{q}}) 
\end{gather}
as $t\to\infty$ uniformly for all $v\in [0,E_0)$ and $\varphi\in\mathbb R$.
\end{Th}

The proof is contained in Section~\ref{Sec3}.

Let us remark that the transformation described in Theorem \ref{Th1} can set to zero some of leading terms in the first equation of \eqref{Veq}. Let $n\geq 1$ be an integer such that
\begin{gather}\label{ass1}
\Lambda_k(v)\equiv 0, \quad k<n, \quad \Lambda_n(v)\not\equiv 0.
\end{gather}

The structure of the transformations \eqref{exch1} and \eqref{exch11} implies that the stability of the trivial solution ${\bf z}(t)\equiv 0$ of system \eqref{FulSys} is associated with the behavior of the trajectories of system \eqref{Veq} with $v(t)$ close to zero, which depends on the structure of the right-hand side of the first equation in system \eqref{Veq}. Note that for non-autonomous systems the stability of solutions depend not only on the leading as $t\to\infty$ and linear terms of equations (see, for example, \cite{OS20arxiv}). With this in mind, consider the following two cases:
\begin{align}
\label{ass21}   &\Lambda_n(v)=\lambda_{n} v\left(1+\mathcal O(v)\right), \quad v\to 0;\\
\label{ass22}     &  
	\begin{cases}
   \Lambda_n(v)=\lambda_{n,m} v^m\left(1+\mathcal O(v)\right), \\
	     \Lambda_{k}(v)=\mathcal O(v^m), \quad n<k<n+l, \\ 
	    \Lambda_{n+l}(v)=\lambda_{n+l} v\left(1+\mathcal O(v)\right), \quad v\to 0.
 \end{cases} 
\end{align}
Here $\lambda_n$, $\lambda_{n,m}$ and $\lambda_{n+l}$ are nonzero constants, $n,m,l$ are integers such that $n\geq 1$, $m\geq 2$, $l\geq 1$. Assumption \eqref{ass21} corresponds to the case, when the leading asymptotic term of the equation as $t\to\infty$ has nonzero linear part. Assumption \eqref{ass22} covers the cases when the leading term is nonlinear in the vicinity of the equilibrium. It follows from the second equation of \eqref{Veq} that $|\varphi(t)|\to\infty$  almost surely as $t\to\infty$ if $v(t)=\mathcal O(1)$. 

Let us specify the definition of stochastic stability that will be used in this paper.
\begin{Def}
The solution ${\bf z}(t)\equiv 0$ of system \eqref{FulSys} is said to be stable
in probability with the weight $\gamma(t)>0$, if there exists $t_0\geq s$ such that 
$\forall\,\varepsilon> 0$ and $\eta > 0$  
$\exists\,\delta > 0$: for all $|{\bf z}_0|< \delta$ the solution ${\bf z}(t)$ of system \eqref{FulSys}
with initial data ${\bf z}(t_0)={\bf z}_0$ satisfies the inequality
\begin{gather}\label{defst}
\mathbb P\left(\sup_{t\geq t_0}\big( |{\bf z}(t)| \gamma(t) \big)>\varepsilon\right)<\eta.
\end{gather}
The solution ${\bf z}(t)\equiv 0$ is said to be unstable in probability with the weight $\gamma(t)$ if it is not stable in probability with this weight.
\end{Def}

This definition modifies slightly classical concept of stability in probability (see, for example, ~\cite[\S 5.3]{RH12}) because of the factor $\gamma(t)$ in the inequality for solutions, which can be considered as the estimate for the norm in the space of continuous functions with the weight $\gamma(t)$. In the case of stability (instability) with the weight $\gamma(t)\equiv 1$, we will say that the solution ${\bf z}(t)\equiv 0$ is stable (unstable) in probability. 

Define the function
\begin{gather*}
\gamma_{n}(t)\equiv 
\begin{cases} 
	\displaystyle \exp\left( \frac{q}{q-n} t^{1-\frac{n}{q}}\right),& n\neq q,\\
	\displaystyle t, & n=q.
\end{cases}
\end{gather*}
It is readily seen that $\gamma_{n}(t)$ grows exponentially (or polynomially) as $t\to\infty$ if $n<q$ (or $n=q$), and $\gamma_{n}(t)$ is bounded if $n>q$. 

Consider first the case \eqref{ass21}.

\begin{Th}\label{Th2}
Let system \eqref{FulSys} satisfy \eqref{zero}, \eqref{lip}, \eqref{H0as}, \eqref{HFBas}, and $n\geq 1$ be an integer such that assumptions \eqref{ass1} and \eqref{ass21} hold. 
\begin{itemize}
	\item 
If $\lambda_n<0$, then for all $\kappa\in (0,1)$ the solution ${\bf z}(t)\equiv 0$ is stable in probability with the weight $(\gamma_{n}(t))^{(1-\kappa)|\lambda_n|/2}$.
\end{itemize}
\end{Th}

From Theorem~\ref{Th2} and the definition of the function $\gamma_{n}(t)$ it follows that the stability is exponential if $n<q$, polynomial if $n=q$ and neutral if $n>q$.

Let us describe the conditions that guarantee the loss of stability. Consider the additional assumption on the intensity of stochastic perturbations in system \eqref{FulSys}:
\begin{gather}\label{ass3}
	\exists\, \mu>0: \quad  | {\hbox{\rm tr}}( {\bf B}^T {\bf B})|\leq \mu^2  t^{-\sigma} |{\bf z}|^2  
\end{gather}
for all $(x,y)\in \mathcal B_r$ and $t\geq  s $ with some constant $\sigma>0$. Then, we have the following.

\begin{Th}\label{Th3}
Let system \eqref{FulSys} satisfy \eqref{zero}, \eqref{lip}, \eqref{H0as}, \eqref{HFBas}, and $1\leq n\leq q$ be an integer such that assumptions \eqref{ass1}, \eqref{ass21} and \eqref{ass3} hold with $\sigma\geq n/q$. 
\begin{itemize}
	\item If $\displaystyle \lambda_n>\delta_{\sigma,n/q}\frac{\mu^2}{2}$, then the solution ${\bf z}(t)\equiv 0$ is unstable in probability.
	\item If $\displaystyle \sigma=\frac{n}{q}$ and $\displaystyle 0< \lambda_n\leq \frac{\mu^2}{2}$, then for all $\kappa>0$ the solution ${\bf z}(t)\equiv 0$ is unstable in probability with the weight 
	$\displaystyle (\gamma_{n}(t))^{({\mu^2}/{2}-\lambda_n+\kappa )/2}$.
\end{itemize}
\end{Th}

Here $\delta_{\sigma,n/q}$ is the Kronecker delta. 

Note that in the case of $\sigma=n/q$ and $0< \lambda_n\leq {\mu^2}/{2}$, Theorem~\ref{Th3} provides only a weak instability of the solution with the weight growing in time. It can be shown that in this case the solution ${\bf z}(t)\equiv 0$ is stable in probability on a finite but asymptotically long time interval as $\mu\to 0$. This property is a variant of a practical stability~\cite{LL61}.

\begin{Th}\label{AsL}
Let system \eqref{FulSys} satisfy \eqref{zero}, \eqref{lip}, \eqref{H0as}, \eqref{HFBas}, and $1\leq n\leq q$ be an integer such that assumptions \eqref{ass1}, \eqref{ass21} and \eqref{ass3} hold with $\sigma=n/q$. 
\begin{itemize}
	\item 
If $\displaystyle 0< \lambda_n\leq \frac{\mu^2}{2}$, then there exists $t_0\geq  s $ such that $\forall\,\varepsilon> 0$ and $\eta > 0$  
$\exists\,\delta > 0$: for all $|{\bf z}_0|< \delta$ the solution ${\bf z}(t)$ of system \eqref{FulSys}
with initial data ${\bf z}(t_0)={\bf z}_0$ satisfies the inequality
\begin{gather}\label{defAsL}
\mathbb P\left(\sup_{0<t-t_0\leq \mathcal T}  |{\bf z}(t)| >\varepsilon\right)<\eta,
\end{gather}
where $\mathcal T= t_0^{n/q}\delta^2(\varepsilon\mu)^{-2}$ if $n<q$, and $\mathcal T=t_0(\exp(\delta^2(\varepsilon\mu)^{-2})-1)$ if $n=q$.
\end{itemize}
\end{Th}

The proofs of Theorems~\ref{Th2}, \ref{Th3} and~\ref{AsL} are contained in Section~\ref{Sec41}. A summary of these results are shown in Table~\ref{Table1}.  
\begin{table}
\begin{tabular}{l|c|c|c}
\hline \rule{0cm}{0.5cm}
  {\bf Assumptions }  &  {\bf Conditions }  & {\bf Stability} & {\bf Ref.}    \\ [0.1cm]
\hline \rule{0cm}{0.5cm}
\eqref{ass1}, \eqref{ass21}, $n<q$ &   & exponentially stable   & Th.~\ref{Th2} \\ [0.1cm]
 \cline{1-1} \cline{3-3}   \rule{0cm}{0.5cm}
\eqref{ass1}, \eqref{ass21}, $n=q$ &  $\lambda_n<0$  & polynomially stable  &  \\ [0.1cm]
 \cline{1-1} \cline{3-3}   \rule{0cm}{0.5cm}
\eqref{ass1}, \eqref{ass21}, $n>q$ &    & stable  &  \\ [0.1cm]
\hline \rule{0cm}{0.5cm}
 \eqref{ass1}, \eqref{ass21}, \eqref{ass3}, $n\leq q$, $\sigma\geq \frac{n}{q}$ & $\lambda_n>\delta_{\sigma,n/q}\frac{\mu^2}{2}$                     & unstable & Th.~\ref{Th3} \\ [0.1cm]
\hline \rule{0cm}{0.5cm}
 \eqref{ass1}, \eqref{ass21}, \eqref{ass3}, $n\leq q$, $\sigma=\frac{n}{q}$ & $0<\lambda_n\leq \frac{\mu^2}{2}$ & practically stable & Th.~\ref{AsL} \\ [0.1cm]
\hline
\end{tabular}
\bigskip
\caption{{\small Stochastic stability of the trivial solution to system \eqref{FulSys} in case \eqref{ass21}.}}
\label{Table1}
\end{table}

Now, consider the case when the leading term in the first equation of \eqref{Veq} is strongly nonlinear as $v\to 0$. Define
\begin{gather}\label{axudef}
\vartheta=\frac{l}{q(m-1)}, \quad 
u_\ast=\left(\frac{|\lambda_{n+l}+\delta_{n+l,q}\vartheta|}{|\lambda_{n,m}|}\right)^{\frac{1}{m-1}}, \quad 
d_\vartheta({\bf z},t;u_\ast):=  t^{\vartheta} H_0(x,y) - u_\ast.
\end{gather}
Then, we have the following.

\begin{Th}\label{Th4}
Let system \eqref{FulSys} satisfy \eqref{zero}, \eqref{lip}, \eqref{H0as}, \eqref{HFBas}, and $n\geq 1$, $l\geq 1$, $m\geq 2$ be integers such that assumptions \eqref{ass1} and \eqref{ass22} hold. 
\begin{itemize}
	\item If $\lambda_{n,m}<0$ and $\lambda_{n+l}<0$, then the solution ${\bf z}(t)\equiv 0$ is stable in probability.
	\item If $n+l\leq q$, $\lambda_{n+l}+\delta_{n+l,q}\vartheta<0$, then the solution ${\bf z}(t)\equiv 0$ is stable in probability with the weight $t^{\vartheta/2}(\gamma_{n+l}(t))^{(1-\kappa)|\lambda_{n+l}+\delta_{n+l,q}\vartheta|/2}$.
	\item If $n+l= q$, $\lambda_{n,m}<0$ and $\lambda_{n+l}>0$, then for all $\varepsilon>0$ and $\eta>0$ there exist $\delta_0>0$ and $t_0>0$ such that the solution ${\bf z}(t)$ of system \eqref{FulSys}
with initial data ${\bf z}(t_0)={\bf z}_0$, $|d_\vartheta({\bf z}_0,t_0;u_\ast)|< \delta_0$, satisfies
\begin{gather}\label{defineq2}
\mathbb P\Big(\sup_{ t\geq t_0 }  |d_\vartheta({\bf z}(t),t;u_\ast)|>\varepsilon\Big)<\eta.
\end{gather}
\end{itemize}
\end{Th}

Note that in the case of \eqref{defineq2} the solution ${\bf z}(t)$ of system \eqref{FulSys}, starting in the vicinity of the equilibrium, has with high probability the following asymptotic behavior: $|{\bf z}(t)|=\mathcal O(t^{-\vartheta/2})$ as $t\to\infty$. This corresponds to a polynomial stability of the equilibrium $(0,0)$.

\begin{Th}\label{Th5}
Let system \eqref{FulSys} satisfy \eqref{zero}, \eqref{lip}, \eqref{H0as}, \eqref{HFBas}, and $n\geq 1$, $l\geq 1$, $m\geq 2$  be integers such that assumptions \eqref{ass1}, \eqref{ass22}, \eqref{ass3} hold with $\sigma\geq (n+l)/q$ and $n+l\leq q$. 
\begin{itemize}
	\item If $\lambda_{n,m}>0$ and $\displaystyle \lambda_{n+l}>\delta_{\sigma,(n+l)/q}\frac{\mu^2}{2}$, then the solution ${\bf z}(t)\equiv 0$ is unstable in probability.
\end{itemize}
\end{Th}

The results obtained for the case when assumption \eqref{ass22} holds are shown in Table~\ref{Table2}. The proofs of the corresponding theorems are contained in Section~\ref{Sec42}.
\begin{table}
\begin{tabular}{l|c|c|c}
\hline \rule{0cm}{0.5cm}
  {\bf Assumptions }  &  {\bf Conditions }  & {\bf Stability} & {\bf Ref.}    \\ [0.1cm]
\hline \rule{0cm}{0.5cm}
\eqref{ass1}, \eqref{ass22} & $\lambda_{n,m}<0$, $\lambda_{n+l}<0$  & stable   & Th.~\ref{Th4} \\ [0.1cm]
 \cline{1-3}  \rule{0cm}{0.5cm}
\eqref{ass1}, \eqref{ass22}, $n+l< q$ &  $\lambda_{n+l}<0$  & exponentially stable  &   \\ [0.1cm]
\cline{1-3} \rule{0cm}{0.5cm}
 \eqref{ass1}, \eqref{ass22}, $n+l= q$ &  $\lambda_{n+l}+\vartheta<0$  & polynomially stable  &   \\ [0.1cm]
\cline{2-2}  \rule{0cm}{0.5cm}
 &  $\lambda_{n,m}<0$, $\lambda_{n+l}>0$  &   &   \\ [0.1cm]
 \hline    \rule{0cm}{0.5cm}
\eqref{ass1}, \eqref{ass22}, \eqref{ass3}, $n+l\leq q$, $\sigma\geq \frac{n+l}{q}$ & $\lambda_{n,m}>0$, $\lambda_{n+l}>\delta_{\sigma,(n+l)/q}\frac{\mu^2}{2}$   & unstable  &  Th.~\ref{Th5} \\ [0.1cm]
\hline
\end{tabular}
\bigskip
\caption{{\small Stochastic stability of the trivial solution to system \eqref{FulSys} in case \eqref{ass22}.}}
\label{Table2}
\end{table}

Finally, consider the case when decaying perturbations lead to the appearance of limit cycles in system \eqref{FulSys}. Let us introduce one more assumption:
\begin{gather}\label{ass30}
\exists\, c\in (0,E_0): \quad \Lambda_n(c)=0, \quad \Lambda'_n(c)\neq 0.
\end{gather}
Then, we have the following.
\begin{Th}\label{Th6}
Let system \eqref{FulSys} satisfy \eqref{zero}, \eqref{lip}, \eqref{H0as}, \eqref{HFBas}, and assumptions \eqref{ass1}, \eqref{ass30} hold with $n=q$. 
\begin{itemize}
	\item If $\Lambda'_n(c)<0$, then for all $\varepsilon>0$ and $\eta>0$ there exist $\delta>0$ and $t_0>0$ such that the solution ${\bf z}(t)$ of system \eqref{FulSys}
with initial data ${\bf z}(t_0)={\bf z}_0$, $|d_0({\bf z}_0,t_0;c)|< \delta$, satisfies
\begin{gather}\label{defineq3}
\mathbb P\Big(\sup_{ t\geq t_0 }  |d_0({\bf z}(t),t;c)|>\varepsilon\Big)<\eta.
\end{gather}
\end{itemize}
\end{Th}
The proof is contained in Section~\ref{Sec43}.

It follows from \eqref{H0as} and \eqref{axudef} that in the case of \eqref{defineq3} the solution of system \eqref{FulSys} has with high probability the following asymptotic behavior: $|{\bf z}(t)|=\sqrt{2c}+ o(1)$ as $t\to\infty$.

\section{Change of variables}\label{Sec3}
In this section, we construct the transformations of variables that reduce system \eqref{FulSys} to \eqref{Veq}.
 
\subsection{Energy-angle variables}
Consider auxiliary $2\pi$-periodic functions $X(\varphi,E)\equiv \hat x(\varphi/\nu(E),E)$ and $Y(\varphi,E)\equiv \hat y(\varphi/\nu(E),E)$. It follows from \eqref{LimSys} that
\begin{gather*}
    \nu(E)\frac{\partial X}{\partial \varphi}=\partial_Y H_0(X,Y), \quad
    \nu(E)\frac{\partial Y}{\partial \varphi}=-\partial_X H_0(X,Y).
\end{gather*}
These functions are used for rewriting system \eqref{FulSys} in the energy-angle variables $(E,\varphi)$. By differentiating the identity $H_0(X(\varphi,E),Y(\varphi,E))\equiv E$ with respect to $E$, we obtain
\begin{gather*}
\det\frac{\partial(X,Y)}{\partial (\varphi,E)}=\begin{vmatrix}
        \partial_\varphi X & \partial_E X\\
        \partial_\varphi Y& \partial_E Y
    \end{vmatrix} = \frac{1}{\nu(E)}\neq 0.
\end{gather*}
Hence, the transformation $(x,y)\mapsto (E,\varphi)$ is invertible for all $E\in [0,E_0]$ and $\varphi\in[0,2\pi)$. Denote by 
\begin{gather}\label{IPhi}
	E=I(x,y),  \quad \varphi=\Phi(x,y)
\end{gather}
the inverse transformation to \eqref{exch1}, and by 
\begin{align}
	\label{Loper}
		\begin{split}
L := &\partial_t  +\partial_y H(x,y,t) \partial_x+\Big(-\partial_x H(x,y,t) + F(x,y,t)\Big) \partial_y   \\
&+\frac{1}{2}
  \Big(\big(B_{1,1}(x,y,t)\big)^2+\big(B_{1,2}(x,y,t)\big)^2\Big)\partial_x^2+\frac{1}{2}\Big(\big(B_{2,1}(x,y,t)\big)^2+\big(B_{2,2}(x,y,t)\big)^2\Big)  \partial_y^2 \\
& + \Big(B_{1,1}(x,y,t)B_{2,1}(x,y,t) +B_{1,2}(x,y,t)B_{2,2}(x,y,t)\Big)\partial_x\partial_y
	\end{split}
\end{align} 
the operator associated with system \eqref{FulSys}. Note that this operator plays the important role in the analysis of stochastic differential equations~\cite[\S 3.3]{RH12}. By applying It\^{o}'s formula to \eqref{IPhi}, it can be shown that, in the new coordinates $(E,\varphi)$, the perturbed system \eqref{FulSys} takes the form
\begin{equation}
	\label{FulSys2}
	\begin{split}
		{dE} =f(E,\varphi,t)dt+\sum_{j=1}^2 \beta_{1,j}(E,\varphi,t) dw_j(t),   \\
		{d\varphi} =\big(\nu(E)+g (E,\varphi,t)\big)dt+\sum_{j=1}^2 \beta_{2,j}(E,\varphi,t) dw_j(t),
	\end{split}
\end{equation}
where
\begin{eqnarray*}
\begin{pmatrix} f \\ g \end{pmatrix}  &\equiv&  -\begin{pmatrix} 0 \\ \nu \end{pmatrix}+
	L \begin{pmatrix}  I \\  \Phi \end{pmatrix}\Big|_{x=X(\varphi,E), y=Y(\varphi,E)},\\
\begin{pmatrix} \beta_{1,j} \\ \beta_{2,j}\end{pmatrix}   & \equiv &    
	\left(B_{1,j} \partial_x +  B_{2,j} \partial_y \right)\begin{pmatrix} I\\ \Phi \end{pmatrix}\Big|_{x=X(\varphi,E), y=Y(\varphi,E)}.
\end{eqnarray*}
It can easily be checked that 
\begin{gather*}
  \partial_x \begin{pmatrix} I \\  \Phi \end{pmatrix} \Big|_{x=X(\varphi,E), y=Y(\varphi,E)} \equiv \begin{pmatrix} -\nu \partial_\varphi Y \\  \nu\partial_E Y \end{pmatrix},\quad \partial_y \begin{pmatrix} I \\  \Phi \end{pmatrix} \Big|_{x=X(\varphi,E), y=Y(\varphi,E)} \equiv \begin{pmatrix}  \nu\partial_\varphi X \\  - \nu\partial_E X \end{pmatrix},\\
  \partial^2_x \begin{pmatrix} I \\  \Phi \end{pmatrix}\Big|_{x=X(\varphi,E), y=Y(\varphi,E)} \equiv \nu (\partial_E Y \partial_\varphi-\partial_\varphi Y\partial_E) \begin{pmatrix} -\nu \partial_\varphi Y \\  \nu\partial_E Y \end{pmatrix},\\ 
\partial^2_y \begin{pmatrix} I \\  \Phi \end{pmatrix}\Big|_{x=X(\varphi,E), y=Y(\varphi,E)} \equiv \nu (\partial_\varphi X \partial_E-\partial_E X\partial_\varphi) \begin{pmatrix} \nu \partial_\varphi X \\  -\nu\partial_E X \end{pmatrix},\\
\partial_x\partial_y \begin{pmatrix} I \\  \Phi \end{pmatrix}\Big|_{x=X(\varphi,E), y=Y(\varphi,E)} \equiv \nu (\partial_\varphi X \partial_E-\partial_E X \partial_\varphi) \begin{pmatrix} -\nu \partial_\varphi Y \\  \nu\partial_E Y \end{pmatrix}.
\end{gather*}
Therefore, the functions $f(E,\varphi,t)$, $g(E,\varphi,t)$ and $\beta_{i,j}(E,\varphi,t)$ are $2\pi$-periodic with respect to $\varphi$. The condition \eqref{zero} implies that $f(0,\varphi,t)\equiv 0$ and $\beta_{1,j}(E,\varphi,t)\equiv 0$. From \eqref{LimSys} and \eqref{H0as} it follows that $X(\varphi,E)=\mathcal O(\sqrt E)$ and $Y(\varphi,E)=\mathcal O(\sqrt E)$ as $E\to 0$ uniformly for all $\varphi\in\mathbb R$. Hence, 
\begin{gather}\label{fgbbes}
f(E,\varphi,t)=\mathcal O(E), \quad  
g(E,\varphi,t)=\mathcal O(1), \quad 
\beta_{1,j}(E,\varphi,t)=\mathcal O(E), \quad 
\beta_{2,j}(E,\varphi,t)=\mathcal O(1)
\end{gather}
 as $E\to 0$ uniformly for all $\varphi\in\mathbb R$ and $t\geq  s $. Moreover, from \eqref{HFBas} it follows that these functions have the following asymptotic expansions: 
\begin{gather*}
f(E,\varphi,t)=\sum_{k=1}^\infty t^{-\frac kq} f_k(E,\varphi), \quad 
g(E,\varphi,t)=\sum_{k=1}^\infty t^{-\frac kq} g_k(E,\varphi), \quad
\beta_{i,j}(E,\varphi,t)=\sum_{k=1}^\infty t^{-\frac kq} \beta_{i,j,k}(E,\varphi)
\end{gather*}
as $t\to\infty$ uniformly for all $E\in [0,E_0]$ and $\varphi\in\mathbb R$, where
\begin{eqnarray*}
\begin{pmatrix} f_{k} \\ g_k \end{pmatrix} &\equiv &
\left(\partial_y H_k \partial_x +
(-\partial_x H_k+F_k)\partial_y \right)\begin{pmatrix} I \\ \Phi \end{pmatrix}\Big|_{x=X(\varphi,E), y=Y(\varphi,E)}\\
	&&+\frac{1}{2}\sum_{i_1+i_2=k} \left\{ (B_{1,1,i_1}B_{1,1,i_2}+B_{1,2,i_1}B_{1,2,i_2})\partial_x^2 
  +  (B_{2,1,i_1}B_{2,1,i_2}+B_{2,2,i_1}B_{2,2,i_2})\partial_y^2 \right. \\ 
	&& +\left.
2(B_{1,1,i_1}B_{2,1,i_2}  +B_{1,2,i_1}B_{2,2,i_2})\partial_x\partial_y  \right \}
\begin{pmatrix} I \\ \Phi \end{pmatrix}\Big|_{x=X(\varphi,E), y=Y(\varphi,E)},\\
\begin{pmatrix} \beta_{1,j,k} \\ \beta_{2,j,k} \end{pmatrix} &\equiv& (B_{1,j,k}\partial_x +B_{2,j,k}\partial_y)\begin{pmatrix} I \\ \Phi \end{pmatrix}\Big|_{x=X(\varphi,E), y=Y(\varphi,E)}.
\end{eqnarray*}

\subsection{Averaging}
It follows from \eqref{FulSys2} that $\varphi(t)$ changes faster than possible variations of $E(t)$ for large values of $t$. The next simplification of the system is associated with averaging the first equation in \eqref{FulSys2} with respect to $\varphi$. This technique is usually used in perturbation theory (see, for example, \cite{BM61,AN84,BDP01,AKN06,DM10}).

Consider a near-identity transformation in the following form:
\begin{gather}\label{VN}
V_N(E,\varphi,t)=E+\sum_{k=1}^N t^{-\frac{k}{q}} v_k(E,\varphi)
\end{gather}
with some integer $N\geq 1$. The coefficients $v_k(E,\varphi)$ are sought in such a way that the drift term of the first equation in transformed system \eqref{Veq} for the new variable $v(t) \equiv  V_N(E(t),\varphi(t),t)$ does not depend explicitly on $\varphi$ at least in the first terms of the asymptotics with some functions $\Lambda_k(v)$ and the remainder satisfying the estimate $\tilde \Lambda_N(v,\varphi,t)=\mathcal O(t^{-{(N+1)}/{q}})$ as $t\to\infty$. Applying It\^{o}'s formula to $V_N(E,\varphi,t)$ along the trajectories of system \eqref{FulSys2} yields
\begin{gather}\label{dVN}
dv=\mathcal L V_N dt+\sum_{j=1}^2\big(\beta_{1,j}\partial_E V_N+\beta_{2,j}\partial_\varphi V_N\big) \, dw_j(t),
\end{gather}
where
\begin{gather}\label{calL}
\mathcal L :=\partial_t + f \partial_E+(\nu+g)\partial_\varphi+\frac{1}{2}(\beta_{1,1}^2+\beta_{1,2}^2)\partial^2_E+(\beta_{1,1}\beta_{2,1}+\beta_{1,2}\beta_{2,2})\partial_E\partial_\varphi+\frac{1}{2}(\beta_{2,1}^2+\beta_{2,2}^2)\partial^2_\varphi
\end{gather}
is the operator associated with system \eqref{FulSys2}. It can easily be checked that 
\begin{equation*}
\begin{split}
\mathcal LV_N=&\sum_{k=1}^\infty t^{-\frac{k}{q}} \Big(f_k+\nu \partial_\varphi v_k-\frac{k-q}{q}v_{k-q}\Big)+\sum_{k=2}^\infty t^{-\frac{k}{q}} \sum_{i_1+i_2=k}(f_{i_1}\partial_E  +g_{i_1}\partial_\varphi )v_{i_2}\\
 & + \frac{1}{2}\sum_{k=3}^\infty t^{-\frac{k}{q}} \sum_{i_1+i_2+i_3=k} \left\{(\beta_{1,1,i_1}\beta_{1,1,i_2}+\beta_{1,2,i_1}\beta_{1,2,i_2})\partial_E^2 \right.\\
& \left. +2(\beta_{1,1,i_1}\beta_{2,1,i_2}+\beta_{1,2,i_1}\beta_{2,2,i_2})\partial_E \partial_\varphi  
+ (\beta_{2,1,i_1}\beta_{2,1,i_2}+\beta_{2,2,i_1}\beta_{2,2,i_2})\partial_\varphi^2  \right\}v_{i_3},
\end{split}
\end{equation*}
where it is assumed that $f_k\equiv g_k\equiv \beta_{i,j,k}\equiv 0$ if $k<1$ and $v_m\equiv 0$ if $m>N$ or $m<1$. Matching \eqref{dVN} with the first equation in \eqref{Veq} gives the following chain of differential equations for determining $v_k(E,\varphi)$:
\begin{gather}\label{vkeq}
\nu(E)\partial_\varphi v_k=\Lambda_k(E)-f_k(E,\varphi)-R_k(E,\varphi), \quad 1\leq k\leq N,
\end{gather}
where each function $R_k$ is expressed through $v_1,\dots,v_{k-1}$. In particular, $R_1\equiv 0$,
\begin{eqnarray*}
R_2&\equiv &(f_1\partial_E+g_1\partial_\varphi)v_1-v_1\partial_E\Lambda_1 -\frac{2-q}{q}v_{2-q},\\
R_3&\equiv & \sum_{i_1+i_2=3}\Big((f_{i_1}\partial_E+g_{i_1}\partial_\varphi)v_{i_2}-v_{i_1}\partial_E\Lambda_{i_2}\Big)-\frac{v_1^2}{2}\partial_E^2 \Lambda_1-\frac{3-q}{q}v_{3-q}\\
&&+\frac{1}{2}\left\{(\beta_{1,1,1}^2+\beta_{1,2,1}^2)\partial_E^2  +2(\beta_{1,1,1}\beta_{2,1,1}+\beta_{1,2,1}\beta_{2,1,1})\partial_E \partial_\varphi  + (\beta_{2,1,1}^2+\beta_{2,2,1}^2)\partial_\varphi^2 \right\} v_{1},\\
\end{eqnarray*}
etc. 
Define 
\begin{gather}\label{lambdak}
\Lambda_k(E)=\langle f_k(E,\varphi)+R_k(E,\varphi)\rangle_\varphi, \quad \langle C(E,\varphi)\rangle_\varphi:=\frac{1}{2\pi}\int\limits_0^{2\pi} C(E,\sigma)\,d\sigma.
\end{gather}
Then, for all $k\in [1,N]$ the right-hand side of \eqref{vkeq} is $2\pi$-periodic in $\varphi$ with zero average. Integrating \eqref{vkeq} yields
\begin{gather*}
v_k(E,\varphi)=-\frac{1}{\nu(E)}\int\limits_0^\varphi \{f_k(E,\sigma)+R_k(E,\sigma)\}_\sigma\, d\sigma+\hat v_k(E),
\end{gather*}
where $\{C(E,\sigma)\}_\sigma=C(E,\sigma)-\langle C(E,\sigma)\rangle_\sigma$ and $\hat v_k(E)$ is chosen such that $\langle v_k(E,\varphi)\rangle_\varphi\equiv 0$. Thus, each function $v_k(E,\varphi)$ is smooth and periodic with respect to $\varphi$. Moreover, it follows from \eqref{vkeq} and \eqref{lambdak} that $v_k(E,\varphi)=\mathcal O(E)$ and $\Lambda_k(E)=\mathcal O(E)$ as $E\to 0$ uniformly for all $\varphi\in\mathbb R$. 

It follows from \eqref{VN} that for all $\epsilon\in (0, 1)$ there exists $t_\ast \geq  s $ such that
\begin{gather}\label{VNineq}
	|V_N(E,\varphi,t)-E|\leq \epsilon E, \quad |\partial_E V_N(E,\varphi,t)-1|\leq \epsilon, \quad |\partial_\varphi V_N(E,\varphi,t)|\leq \epsilon E
\end{gather}
for all $E\in [0,E_0]$, $\varphi\in \mathbb R$ and $t\geq t_\ast$. Hence, the mapping $(E,\varphi,t)\mapsto (v,\varphi,t)$ is invertible for all $v\in [0,v_\ast]$, $\varphi\in [0,2\pi)$ and $t\geq t_\ast$ with $v_\ast=(1-\epsilon) E_0$. 

Let $E=\mathcal E(v,\varphi,t)$ be the inverse transformation to \eqref{VN}.
Then,
\begin{equation*}
\begin{split}
	& \alpha_{1,j}(v,\varphi,t)\equiv \left(\beta_{1,j}(E,\varphi,t)\partial_E  +\beta_{2,j}(E,\varphi,t)\partial_\varphi \right)V_N(E,\varphi,t)\Big|_{E=\mathcal E(v,\varphi,t)},\\
	& \tilde \Lambda_N(v,\varphi,t) \equiv -\sum_{k=1}^N t^{-\frac{k}{q}} \Lambda_k(v)+ \mathcal L V_N(E,\varphi,t)\Big|_{E=\mathcal E(v,\varphi,t)}.
\end{split}
\end{equation*}
The second equation in \eqref{FulSys2} does not change significantly under the transformation \eqref{VN}. In particular, $\tilde G_N(v,\varphi,t)\equiv \nu(\mathcal E(v,\varphi,t))-\nu(v)+g(\mathcal E(v,\varphi,t),\varphi,t)$ and $\alpha_{2,j}(v,\varphi,t)\equiv \beta_{2,j}(\mathcal E(v,\varphi,t),\varphi,t)$ in the second equation of the transformed system \eqref{Veq}. It can easily be checked that
\begin{gather*}
\tilde \Lambda_N(v,\varphi,t)=\mathcal O(t^{-\frac{N+1}{q}}), \quad \tilde G_N(v,\varphi,t)=\mathcal O(t^{-\frac{1}{q}}), \quad \alpha_{i,j}(v,\varphi,t)=\mathcal O(t^{-\frac{1}{q}})
\end{gather*}
as $t\to\infty$ uniformly for all $v\in [0,v_\ast]$ and $\varphi\in\mathbb R$. From \eqref{fgbbes} it follows that
\begin{gather*}
\tilde \Lambda_{N}(v,\varphi,t)=\mathcal O(v), \quad
\tilde G_{N}(v,\varphi,t)=\mathcal O(1), \quad 
\alpha_{1,j}(v,\varphi,t)=\mathcal O(v), \quad 
\alpha_{2,j}(v,\varphi,t)=\mathcal O(1)
\end{gather*} 
as $v\to 0$ uniformly for all $t\geq t_\ast$ and $\varphi\in\mathbb R$.

Thus, we obtain the proof of Theorem~\ref{Th1}.

\section{Stability analysis in case \eqref{ass21}}
\label{Sec41}

\subsection*{Proof of Theorem~\ref{Th2}}
Let us fix $\kappa\in (0,1)$ and consider 
\begin{gather*}
	U_{n,\kappa}(x,y,t)= \gamma_{n,\kappa}(t) V_N(I(x,y),\Phi(x,y),t)
\end{gather*}
with $N=n$ and the positive function $\gamma_{n,\kappa}(t)\equiv (\gamma_{n}(t))^{(1-\kappa)|\lambda_n|}$ as a Lyapunov function candidate for system \eqref{FulSys}. The function $V_N(E,\varphi,t)$ having the form \eqref{VN} was constructed in the previous section. The functions $I(x,y)$ and $\Phi(x,y)$ are defined by \eqref{IPhi}. It follows from \eqref{H0as} and \eqref{VNineq} that for all $\epsilon>0$ there exists $r_\ast\in (0,r]$ such that 
\begin{gather}\label{Uineq}
(1-\epsilon)^2 \frac{|{\bf z}|^2}{2} \leq V_n(I(x,y),\Phi(x,y),t)\leq (1+\epsilon)^2 \frac{|{\bf z}|^2}{2}
\end{gather} 
for all $t\geq t_\ast$ and $(x,y)\in \mathcal B_{r_\ast}\subseteq \mathcal D(E_0)$. 
It can easily be checked that 
\begin{gather*}
LU_{n,\kappa}(x,y,t)\equiv \frac{\gamma_{n,\kappa}'(t)}{\gamma_{n,k}(t)}U_{n,\kappa}(x,y,t)+\gamma_{n,\kappa}(t) \mathcal L V_n(E,\varphi,t)\Big|_{E=I(x,y), \varphi=\Phi(x,y)}.
\end{gather*}
Recall that $L$ and $\mathcal L$ are defined by \eqref{Loper} and \eqref{calL}, respectively.
Note that $\gamma_{n,\kappa}'(t)/\gamma_{n,\kappa}(t)\equiv (1-\kappa)|\lambda_n|t^{-n/q}$. 
Combining this with \eqref{ass1}, we obtain
\begin{gather*}
L U_{n,\kappa} (x,y,t) = t^{-\frac{n}{q}} U_{n,\kappa}(x,y,t) \Big (-\kappa |\lambda_n|+\mathcal O(|{\bf z}|^2)+\mathcal O(t^{-\frac{1}{q}})\Big)
\end{gather*}
as $t\to\infty$ and ${\bf z}\to 0$. Hence, for all $\epsilon\in (0,1)$ there exist $t_0\geq t_\ast$ and $r_0\leq r_\ast$ such that
\begin{gather}\label{LUineq}
LU_{n,\kappa}(x,y,t)\leq -t^{-\frac{n}{q}} U_{n,\kappa}(x,y,t) (1-\epsilon) \kappa |\lambda_n| \leq -t^{-\frac{n}{q}} \gamma_{n,\kappa}(t) \frac{|{\bf z}|^2}{2} (1-\epsilon)^3  \kappa |\lambda_n| \leq 0
\end{gather}
for all $t\geq t_0$ and $(x,y)\in \mathcal B_{r_0}$.

Fix the parameters $0<\varepsilon<r_0\sqrt{\gamma_{n,\kappa}(t_0)}$ and $0<\eta<1$. Let ${\bf z}(t)\equiv (x(t),y(t))^T$ be a solution of system \eqref{FulSys} with initial data $(x(t_0),y(t_0))\in \mathcal B_\delta$, $0<\delta<\varepsilon$, and $\tau_\varepsilon$ be the first exit time of $\tilde{\bf z}(t)$ from the domain $\mathcal B_\varepsilon$ as $t>t_0$, where $\tilde{\bf z}(t)\equiv {\bf z}(t)\sqrt{\gamma_{n,\kappa}(t)}$.
Define the function $\tau_\varepsilon(t)=\min\{\tau_\varepsilon,t\}$, then $\tilde{\bf z}(\tau_\varepsilon(t))$ is the process stopped at a first exit time from $\mathcal B_\varepsilon$. From \eqref{LUineq} it follows that $U(x(\tau_\varepsilon(t)),y(\tau_\varepsilon(t)),\tau_\varepsilon(t))$ is a non-negative supermartingale (see, for example,~\cite[\S 5.2]{RH12}). Hence, using \eqref{Uineq} and the definition of $\tau_\varepsilon(t)$, we get the following:
\begin{align*}
\mathbb P\left(\sup_{t\geq t_0} |\tilde {\bf z}(t)|\geq \varepsilon\right)
& =	\mathbb P\left(\sup_{t\geq t_0}\big(|{\bf z}(\tau_\varepsilon(t))|^2\gamma_{n,\kappa}(\tau_\varepsilon(t))\big)\geq \varepsilon^2\right)\\
& \leq \mathbb P\left(\sup_{t\geq t_0}U_{n,\kappa}\big(x(\tau_\varepsilon(t)),y(\tau_\varepsilon(t)),\tau_\varepsilon(t)\big)\geq (1-\epsilon)^2\frac{\varepsilon^2}{2}\right)\\
& \leq \frac{ U_{n,\kappa}(x(t_0),y(t_0),t_0)}{(1-\epsilon)^2 \varepsilon^2/2}.
\end{align*}
The last estimate follows from Doob's supermartingale inequality. From \eqref{Uineq} it follows that 
\begin{gather*}
	U_{n,\kappa}(x(t_0),y(t_0),t_0)\leq \gamma_{n,\kappa}(t_0)(1+\epsilon)^2 \frac{\delta^2}{2}.
\end{gather*}
Hence, taking $\delta=\varepsilon \sqrt{\eta/\gamma_{n,\kappa}(t_0)} (1-\epsilon)/(1+\epsilon)$, we obtain \eqref{defst} with $\gamma(t)\equiv \sqrt{\gamma_{n,\kappa}(t)}$. 
\qed

\subsection*{Proof of Theorem~\ref{Th3}}
The proof of instability is based on constructing suitable Lyapunov function for system \eqref{FulSys} and reproduces the arguments from \cite[\S 5.4]{RH12} with some modifications.

1. Let $\lambda_n>\delta_{\sigma,n/q}(\mu^2/2)$. Consider
\begin{gather}
\label{UnFun}
U_1(x,y,t)= \log \left((1+\epsilon)^2\frac{ r^2}{2}\right)-\log V_N(I(x,y),\Phi(x,y),t)
\end{gather}
with $N=n$ as a Lyapunov function candidate for system \eqref{FulSys}. From \eqref{Uineq} it follows that $U_1(x,y,t)\geq \log (r/|{\bf z}|)^2\geq 0$ for all $t\geq t_\ast$ and $(x,y)\in\mathcal B_{r_\ast}$. It can easily be checked that 
\begin{align*}
	LU_1(x,y,t)\equiv & 
	\frac{1}{2}{\hbox{\rm  tr}} \big({\bf B}^T {\bf M} {\bf B}\big) - \frac{\mathcal L V_n(E,\varphi,t)}{V_n(E,\varphi,t)}\Big|_{E=I(x,y), \varphi=\Phi(x,y)},
\end{align*}
where
\begin{gather*}\left.
{\bf M}({\bf z},t)\equiv \big[V_n(I(x,y),\Phi(x,y),t)\big]^{-2}
\begin{pmatrix}
(\partial_xV_n)^2 &  \partial_xV_n \partial_yV_n \\
 \partial_xV_n \partial_yV_n  & (\partial_yV_n)^2
\end{pmatrix}\right|_{E=I(x,y), \varphi=\Phi(x,y)}.
\end{gather*}
It follows from \eqref{H0as}, \eqref{ass3}, \eqref{VNineq} and \eqref{Uineq} that for all $\epsilon\in (0,1)$ there exist $r_1\in (0,r_\ast]$ and $t_1\geq t_\ast$ such that
\begin{gather*}
\big|{\hbox{\rm tr}} ({\bf B}^T {\bf M} {\bf B})\big |\leq  \mu^2 t^{-\sigma} \left(\frac{1+\epsilon}{1-\epsilon}\right)^4
\end{gather*}
for all $t\geq t_1$ and $(x,y)\in \mathcal B_{r_1}$. Combining this with \eqref{ass21}, we obtain 
\begin{gather}\label{LUineq1}
LU_1(x,y,t)\leq t^{-\frac{n}{q}} \left(-\lambda_n+\frac{\mu^2}{2}t^{-(\sigma-\frac{n}{q})}\left(\frac{1+\epsilon}{1-\epsilon}\right)^4 + \mathcal O(|{\bf z}|^2)+\mathcal O\big(t^{-\frac{1}{q}}\big)\right)
\end{gather}
as $t\to\infty$ and ${\bf z}\to 0$. Choosing $\epsilon>0$ small enough ensures that there exist $t_2\geq t_1$ and $r_2\leq r_1$ such that $LU_1(x,y,t)\leq 0$ for all $t\geq t_2$ and $(x,y)\in \mathcal B_{r_2}$. 

Consider also the auxiliary function
\begin{gather}\label{AuxFun}
U_\ast(x,y,t)= (1+\epsilon)^2\frac{r^2}{2}-V_N(I(x,y),\Phi(x,y),t)
\end{gather}
with $N=n$.
It is readily seen that 
\begin{align*}
LU_\ast(x,y,t)&  =- t^{-\frac{n}{q}}V_n\big(I(x,y),\Phi(x,y),t\big) \Big(\lambda_n +\mathcal O(|{\bf z}|^2)+\mathcal O\big(t^{-\frac{1}{q}}\big)\Big)
\end{align*}
as $t\to\infty$ and ${\bf z} \to 0$. Hence, for all $\epsilon\in (0,1)$ there exist $r_3\in (0,r_\ast]$ and $t_3\geq t_\ast$ such that 
\begin{gather}\label{Wineq}
U_\ast(x,y,t)\geq 0, \quad LU_\ast(x,y,t)\leq -t^{-\frac{n}{q}} \frac{|{\bf z}|^2}{2} (1-\epsilon)^3\lambda_n \leq 0
\end{gather}
for all $t\geq t_3$ and $|{\bf z}|\leq r_3$.

Define $r_0=\min\{r_2,r_3\}$, $t_0=\max\{t_2,t_3\}$ and fix $\varepsilon\in (0,r_0]$.  Let ${\bf z}(t)\equiv (x(t),y(t))^T$ be a solution of system \eqref{FulSys} with initial data $(x(t_0),y(t_0))\in \mathcal B_\varepsilon$ and let $\delta\in (0,|{\bf z}_0|)$. Define
\begin{gather*}
  \tau_{\delta,\varepsilon}=\inf\left\{t>t_0:  |  {\bf z}(t)|=\delta \ \ \text{or} \ \  | {\bf z}(t)|=\varepsilon\right\},\\ 
   \tau_{\delta}=\inf\left\{t> t_0:  | {\bf z}(t)|=\delta\right\}, \\ 
   \tau_{\varepsilon}=\inf\left\{t> t_0:   |  {\bf z}(t)| =\varepsilon\right\}.
\end{gather*}
We see that $\tau_{\delta,\varepsilon}$, $\tau_{\varepsilon}$ and $\tau_{\delta}$ are the first times of reaching the solution ${\bf z}(t)$ to the boundaries of the domains $\overline{\mathcal B_\varepsilon\setminus \mathcal B_\delta}$, $\mathcal B_\varepsilon$ and $\mathcal B_\delta$, respectively. Define $\tau_{\delta,\varepsilon}(t)=\min\{\tau_{\delta,\varepsilon},t\}$. From the properties of the function $U_1(x,y,t)$ it follows that $U_1(x(\tau_{\delta,\varepsilon}(t)),y(\tau_{\delta,\varepsilon}(t)),\tau_{\delta,\varepsilon}(t))$ is a non-negative supermartingale and 
\begin{gather}\label{expest}
\mathbb {E}\left[U_1\big(x(\tau_{\delta,\varepsilon}(t)),y(\tau_{\delta,\varepsilon}(t)),\tau_{\delta,\varepsilon}(t)\big)\right]\leq U_1(x(t_0),y(t_0),t_0)
\end{gather}
for all $t\geq t_0$. Note that the constructed function $U_\ast(x,y,t)$, satisfying estimates \eqref{Wineq} with $n/q\leq 1$, guarantees the recurrence of ${\bf z}(t)$ (see~\cite[Theorem 3.9]{RH12}). This implies that $\mathbb{P}(\tau_{\delta,\rho}<\infty)=1$. Hence, by letting $t\to\infty$ in \eqref{expest}, we get 
\begin{gather}\label{chineq}
\begin{split}
U_1(x(t_0),y(t_0),t_0)&\geq 
	\mathbb {E}\left[U_1\big(x(\tau_{\delta,\varepsilon}),y(\tau_{\delta,\varepsilon}),\tau_{\delta,\varepsilon}\big)\right]\\
	&\geq \mathbb {E}\left[U_1\big(x(\tau_{\delta }),y(\tau_{\delta }),\tau_{\delta }\big) \chi_{\{\tau_\delta<\tau_\varepsilon\}}\right]\\
	&\geq \mathbb P(\tau_\delta<\tau_\varepsilon) \inf_{t\geq t_0, \, |{\bf z}|=\delta} U_1(x,y,t) \geq \mathbb P\left(\sup_{t_0\leq t\leq \tau_\delta} |{\bf z}(t) |<\varepsilon\right) \log \left(\frac{r}{\delta}\right)^2,
\end{split}
\end{gather}
where $\chi_{\{\tau_\delta<\tau_\varepsilon\}}$ is the indicator function of the set $\{\omega\in\Omega: \tau_\delta<\tau_\varepsilon\}$. Note that the point $(0,0)^T$ is inaccessible to sample paths of the process  ${\bf z}(t)$ (see~\cite[Lemma 5.3]{RH12}). Hence, $\tau_\delta\to\infty$ as $\delta\to 0$ with probability one. Combining this with \eqref{chineq}, we obtain 
\begin{gather*}
\mathbb P\left(\sup_{t\geq t_0} |{\bf z}(t) |<\varepsilon\right) = 0.
\end{gather*}
This implies that the solution ${\bf z}(t)\equiv 0$ is not stable in probability with the weight $\gamma(t)\equiv 1$.

2. Now, let $\sigma=n/q$ and $0< \lambda_n\leq \frac{\mu^2}{2}$. Fix the parameter $\kappa>0$ and consider
\begin{gather*}
U_2(x,y,t)= U_1(x,y,t)-\log  \gamma_{n,\kappa}(t)+K
\end{gather*}
with some parameter $K$ and the positive function $\gamma_{n,\kappa}(t)\equiv (\gamma_{n}(t))^{\mu^2/2-\lambda_n+\kappa}$ as a Lyapunov function candidate for system \eqref{FulSys}. Note that $\big(\log\gamma_{n,\kappa}(t)\big)'\equiv  (\mu^2/2-\lambda_n+\kappa) t^{-n/q}$. 
Combining this with \eqref{LUineq1}, we obtain
\begin{gather*}
	LU_2(x,y,t)\leq t^{-\frac{n}{q}} \left(-\kappa+\frac{\mu^2}{2}\left[\left(\frac{1+\epsilon}{1-\epsilon}\right)^4 -1\right]+ \mathcal O(|{\bf z}|^2)+\mathcal O\big(t^{-\frac{1}{q}}\big)\right)
\end{gather*}
as $t\to\infty$ and ${\bf z}\to 0$. Choosing $\epsilon>0$ small enough ensures that there exist $t_\kappa\geq t_1$ and $r_\kappa\leq r_1$ such that $LU_2(x,y,t)\leq 0$ for all $t\geq t_\kappa$ and $(x,y)\in \mathcal B_{r_\kappa}$.

Define $r_0=\min\{r_\kappa,r_3\}\leq r$, $t_0=\max\{t_\kappa,t_3\}$ and fix $0<\varepsilon< r_0 \sqrt{\gamma_{n,\kappa}(t_0)}$.  Let ${\bf z}(t)\equiv (x(t),y(t))^T$ be a solution of system \eqref{FulSys} with initial data ${\bf z}_0=(x(t_0),y(t_0))^T$ such that $ |{\bf z}_0|\leq \varepsilon/\sqrt{\gamma_{n,\kappa}(t_0)}$ and let $0<\delta< \sqrt{\gamma_{n,\kappa}(t_0)} |{\bf z}_0|$. 
Define 
\begin{gather*}
 \tilde\tau_{\delta,\varepsilon}=\inf\left\{t> t_0:  |\tilde {\bf z}(t)|=\delta \ \ \text{or} \ \  |\tilde {\bf z}(t)|=\varepsilon\right\},\\
 \tilde\tau_{\delta}=\inf\left\{t> t_0:  |\tilde {\bf z}(t)|=\delta\right\}, \\ 
 \tilde\tau_{\varepsilon}=\inf\left\{t> t_0:   |\tilde {\bf z}(t)| =\varepsilon\right\}
\end{gather*}
and $\tilde \tau_{\delta,\varepsilon}(t)=\min\{\tilde{\tau}_{\delta,\varepsilon},t\}$, where $\tilde{\bf z}(t)\equiv {\bf z}(t)\sqrt{\gamma_{n,\kappa}(t)}$. Choose $K=\log \gamma_{n,\kappa}(t_0)$, then 
\begin{gather*}
	U_2\big(x(\tilde \tau_{\delta,\varepsilon}(t)),y(\tilde \tau_{\delta,\varepsilon}(t)),\tilde \tau_{\delta,\varepsilon}(t)\big)\geq \log \left(\frac{\gamma_{n,\kappa}(t_0) r^2}{\varepsilon^2}\right)\geq 0.
\end{gather*}
Hence, $U_2(x(\tilde \tau_{\delta,\varepsilon}(t)),y(\tilde \tau_{\delta,\varepsilon}(t)),\tilde \tau_{\delta,\varepsilon}(t))$ is a non-negative supermartingale and 
\begin{gather}\label{expest2}
\mathbb {E}\left[U_2\big(x(\tilde \tau_{\delta,\varepsilon}(t)),y(\tilde \tau_{\delta,\varepsilon}(t)),\tilde \tau_{\delta,\varepsilon}(t)\big)\right]\leq U_2(x(t_0),y(t_0),t_0)
\end{gather}
for all $t\geq t_0$. Arguing as above, we obtain $\mathbb{P}(\tilde \tau_{\delta,\varepsilon}<\infty)=1$, and by letting $t\to\infty$ in \eqref{expest2}, we get 
\begin{gather}\label{chineq2}
\begin{split}
U_2(x(t_0),y(t_0),t_0)&\geq 
	\mathbb {E}\left[U_2\big(x(\tilde{\tau}_{\delta,\varepsilon}),y(\tilde{\tau}_{\delta,\varepsilon}),\tilde{\tau}_{\delta,\varepsilon}\big)\right]\\
	&\geq \mathbb {E}\left[U_2\big(x(\tilde{\tau}_{\delta }),y(\tilde{\tau}_{\delta }),\tilde{\tau}_{\delta }\big) \chi_{\{\tilde{\tau}_\delta<\tilde{\tau}_\varepsilon\}}\right]\\
	&\geq  \log \left(\frac{ \gamma_{n,\kappa}(t_0) r^2}{\delta^2 }\right) \mathbb P(\tilde{\tau}_\delta<\tilde{\tau}_\varepsilon)  
		=  \log \left(\frac{\gamma_{n,\kappa}(t_0) r^2}{\delta^2 }\right) \mathbb P\left(\sup_{t_0\leq t\leq \tilde\tau_\delta} |\tilde{\bf z}(t) |<\varepsilon\right).
\end{split}
\end{gather}
Since $\tilde \tau_\delta\to\infty$ as $\delta\to 0$ with probability one, it follows from \eqref{chineq2}
that 
\begin{gather*}
\mathbb P\left(\sup_{t\geq t_0} |\tilde{\bf z}(t) |<\varepsilon\right) =0.
\end{gather*}
Hence, the solution ${\bf z}(t)\equiv 0$ is unstable in probability with the weight $\gamma(t)\equiv \sqrt{\gamma_{n,\kappa}(t)}$.
\qed

\subsection*{Proof of Theorem~\ref{AsL}}
Consider the function $U(x,y,t)= V_N(I(x,y),\Phi(x,y),t)$ with $N=n$. It follows from \eqref{RemEst}, \eqref{ass21} and \eqref{Uineq} that for all $\epsilon\in(0,1)$ there exist $t_0\geq t_\ast$ and $r_0\leq r_\ast$ such that 
\begin{gather*}
(1-\epsilon)^2\frac{|{\bf z}|^2}{2} \leq U(x,y,t)\leq (1+\epsilon)^2\frac{|{\bf z}|^2}{2},\quad
LU(x,y,t)\leq t^{-\frac{n}{q}} (1+\epsilon)^2\frac{\mu^2 |{\bf z}|^2}{2}
\end{gather*}
for all $t\geq t_0$ and $(x,y)\in \mathcal B_{r_0}$. 

Fix the parameters $\varepsilon\in (0,r_0)$ and $\eta\in (0,1)$. Let ${\bf z}(t)\equiv (x(t),y(t))^T$ be a solution of system \eqref{FulSys} with initial data $(x(t_0),y(t_0))\in\mathcal B_\delta$ with $\delta\in (0,\varepsilon)$, and $\tau_D$ be the first exit time of the process $(x(t),y(t),t)$ from the domain 
\begin{gather*}
D(\varepsilon,t_0,\mathcal T)=\{(x,y,t): (x,y)\in\mathcal B_\varepsilon, \ \ t_0\leq t\leq t_0+\mathcal T\}.
\end{gather*}
The parameters $\delta>0$ and $\mathcal T>0$ will be determined below. Define the function $\tau_D(t)=\min\{\tau_D,t\}$. 

1. Let $n<q$. Consider
\begin{gather*}
	U_1(x,y,t)=  U(x,y,t)+t_0^{-\frac{n}{q}}\frac{\mu^2\varepsilon^2}{2}(1+\epsilon)^2 (\mathcal T+t_0-t)
\end{gather*}
as a Lyapunov function candidate for system \eqref{FulSys}. It can easily be checked that
\begin{gather*}
U_1(x,y,t)\geq U(x,y,t)\geq 0, \quad
LU_1(x,y,t)\leq t_0^{-\frac{n}{q}}(1+\epsilon)^2 \frac{\mu^2 (|{\bf z}|^2-\varepsilon^2)}{2}\leq 0
\end{gather*}
for all $(x,y,t)\in D(\varepsilon,t_0,\mathcal T)$. Hence $U_1(x(\tau_D(t)),y(\tau_D(t)),\tau_D(t))$ is a non-negative supermartingale and the following estimates hold:
\begin{align*}
\mathbb P\left(\sup_{t_0\leq t\leq t_0+\mathcal T} |{\bf z}(t)| \geq \varepsilon\right)
& =	\mathbb P\left(\sup_{t\geq t_0} |{\bf z}(\tau_D(t))|^2 \geq \varepsilon^2\right)\\
& \leq \mathbb P\left(\sup_{t\geq t_0}U_1(x(\tau_D(t)),y(\tau_D(t)),\tau_D(t))\geq (1-\epsilon)^2\frac{\varepsilon^2}{2}\right)  \leq \frac{ U_1(x(t_0),y(t_0),t_0)}{(1-\epsilon)^2 \varepsilon^2/2}.
\end{align*}
Since
$ U_1(x(t_0),y(t_0),t_0)\leq  (1+\epsilon)^2(\delta^2+\mu^2\varepsilon^2t_0^{-{n}/{q}}\mathcal T)/2$, we can take
 $\delta=\varepsilon \sqrt{C\eta}$ and $\mathcal T= \mu^{-2} t_0^{n/q} C\eta = \mu^{-2} t_0^{n/q}(\delta/\varepsilon)^2 $ with
$C= ((1-\epsilon)/(1+\epsilon))^2/2$ to obtain \eqref{defAsL}.

2. Let $n=q$. Using 
\begin{gather*}
	U_2(x,y,t)=  U(x,y,t)+\frac{\mu^2\varepsilon^2}{2} (1+\epsilon)^2 \log\left(\frac{\mathcal T+t_0}{t}\right)
\end{gather*}
as a Lyapunov function candidate, we obtain 
\begin{gather*}
U_2(x,y,t)\geq U(x,y,t)\geq 0, \quad
LU_2(x,y,t)\leq t^{-1}(1+\epsilon)^2 \frac{\mu^2 (|{\bf z}|^2-\varepsilon^2)}{2}\leq 0
\end{gather*}
for all $(x,y,t)\in D(\varepsilon,t_0,\mathcal T)$. This implies that $U_2(x(\tau_D(t)),y(\tau_D(t)),\tau_D(t))$ is a non-negative supermartingale and, as in the previous case, we obtain the estimate
\begin{align*}
\mathbb P\left(\sup_{t_0\leq t\leq t_0+\mathcal T} |{\bf z}(t)| \geq \varepsilon\right)
 \leq \frac{ U_2(x(t_0),y(t_0),t_0)}{(1-\epsilon)^2 \varepsilon^2/2}.
\end{align*}
It is readily seen that 
$ U_2(x(t_0),y(t_0),t_0)\leq  (1+\epsilon)^2(\delta^2+\mu^2\varepsilon^2\log(1+\mathcal T/t_0))/2$.
Taking $\delta=\varepsilon \sqrt{C\eta }$ and $\mathcal T= t_0(\exp(C \mu^{-2}\eta)-1) $, we obtain \eqref{defAsL}.
\qed

\section{Stability analysis in case \eqref{ass22}}
\label{Sec42}
\subsection*{Proof of Theorem~\ref{Th4}}
1. Let $\lambda_{n,m}<0$ and $\lambda_{n+l}<0$. Using
\begin{gather*}
	U_1(x,y,t)=  V_N(I(x,y),\Phi(x,y),t)
\end{gather*}
with $N=n+l$ as a Lyapunov function candidate for system \eqref{FulSys}, we obtain
\begin{gather}\label{LU1}
	\begin{split}
	LU_1(x,y,t)\equiv & \mathcal L V_{n+l}(E,\varphi,t)\Big|_{E=I(x,y), \varphi=\Phi(x,y)}\\
	 = &t^{-\frac{n}{q}} \big(U_1(x,y,t)\big)^m \left( \lambda_{n,m}+\mathcal O(U_1)+\mathcal O(t^{-\frac{1}{q}}) \right)\\
	&+t^{-\frac{n+l}{q}} U_1(x,y,t)\left(\lambda_{n+l}+\mathcal O (U_1)+\mathcal O(t^{-\frac{1}{q}})\right)
\end{split}
\end{gather}
as $t\to\infty$ and $|{\bf z}|\to 0$. 
Therefore, for all $\epsilon>0$ there exists $r_0\in (0,r]$ and $t_0\geq t_\ast$ such that 
\begin{gather}\label{Uineq0}
(1-\epsilon)^2\frac{|{\bf z}|^2}{2} \leq U_1(x,y,t)\leq (1+\epsilon)^2\frac{|{\bf z}|^2}{2}, \quad LU_1(x,y,t)\leq 0
\end{gather} 
for all $t\geq t_0$ and $(x,y)\in \mathcal B_{r_0}\subseteq \mathcal D(E_0)$.

Fix the parameters $\varepsilon\in (0,r_0)$ and $\eta\in (0,1)$. Let ${\bf z}(t)\equiv (x(t),y(t))^T$ be a solution of system \eqref{FulSys} with initial data $(x(t_0),y(t_0))\in\mathcal B_\delta$ with $\delta\in (0,\varepsilon)$, and $\tau_\varepsilon$ be the first exit time of ${\bf z}(t)$ from the domain $\mathcal B_\varepsilon$ as $t>t_0$.
Define the function $\tau_\varepsilon(t)=\min\{\tau_\varepsilon,t\}$. From \eqref{Uineq0} it follows that $U_1(x(\tau_\varepsilon(t)),y(\tau_\varepsilon(t)),\tau_\varepsilon(t))$ is a non-negative supermartingale, and the following estimates hold:
\begin{align*}
\mathbb P\left(\sup_{t\geq t_0} |{\bf z}(t)| \geq \varepsilon\right)
& =	\mathbb P\left(\sup_{t\geq t_0} |{\bf z}(\tau_\varepsilon(t))|^2 \geq \varepsilon^2\right)\\
& \leq \mathbb P\left(\sup_{t\geq t_0}U_1\big(x(\tau_\varepsilon(t)),y(\tau_\varepsilon(t)),\tau_\varepsilon(t)\big)\geq (1-\epsilon)^2\frac{\varepsilon^2}{2}\right) \leq \frac{ U_1(x(t_0),y(t_0),t_0)}{(1-\epsilon)^2 \varepsilon^2/2}.
\end{align*}
Since $U_1(x(t_0),y(t_0),t_0)\leq (1+\epsilon)^2 \delta^2/2$, we can take $\delta=\varepsilon \sqrt{\eta} (1-\epsilon)/(1+\epsilon)$ to obtain \eqref{defst} with $\gamma(t)\equiv 1$. 

2. Let $n+l\leq q$ and $\lambda^\ast_{n+l}:=\lambda_{n+l}+\delta_{n+l,q}\vartheta<0$  with $\vartheta={l/(q(m-1))}>0$. Consider 
\begin{gather*}
	U_2(x,y,t)= \gamma_{n+l,\kappa,\vartheta}(t) U_1(x,y,t)
\end{gather*}
with the positive function 
$\gamma_{n+l,\kappa,\vartheta}(t)\equiv t^{\vartheta}(\gamma_{n+l}(t))^{(1-\kappa)|\lambda^\ast_{n+l}|}$  as a Lyapunov function candidate. Note that 
\begin{gather*}
	\big(\log\gamma_{n+l,\kappa,\vartheta}(t)\big)'\equiv \vartheta t^{-1}+ (1-\kappa)|\lambda^\ast_{n+l}|t^{-\frac{n+l}{q}}.
\end{gather*}
Combining this with \eqref{LU1}, we obtain
\begin{gather*}
\begin{split}
	LU_2(x,y,t)\equiv & \big(\log\gamma_{n+l,\kappa,\vartheta}(t)\big)' U_2(x,y,t)+\gamma_{n+l,\kappa,\vartheta}(t) L U_1(x,y,t)\\
	 = &t^{-\frac{n+l}{q}}  U_2(x,y,t)\left( -\kappa |\lambda^\ast_{n+l}|+\mathcal O(t^{-\frac{\varkappa}{q}}) \right)
	\end{split}
\end{gather*}
 as $t\to\infty$ and ${\bf z}\to 0$, where $\varkappa=\min\{1,q\vartheta,q(1-\kappa)|\lambda^\ast_{n+l}|\}$. 
Hence, for all $\epsilon>0$ there exists $t_0\geq t_\ast$ and $r_0\leq r_\ast$ such that 
\begin{gather*}
LU_2(x,y,t)\leq -t^{-\frac{n+l}{q}}(1-\epsilon)U_2(x,y,t) \kappa|\lambda^\ast_{n+l}|\leq 0
\end{gather*}
for all $t\geq t_0$ and $(x,y)\in \mathcal B_{r_0}\subseteq \mathcal D(E_0)$. Repeating the proof of Theorem~\ref{Th2}, it can be shown that for all $\varepsilon>0$ and $\eta>0$ there exists $\delta>0$ such that the equilibrium $(0,0)$ of system \eqref{FulSys} is stable in probability with the weight $\gamma(t)\equiv \sqrt{\gamma_{n+l,\kappa,\vartheta}(t)}$. 

3. Finally, consider the case $n+l=q$, $\lambda_{n+l}>0$ and $\lambda_{n,m}<0$. Fix the parameters $\varepsilon>0$ and $\eta>0$. Define the domain
\begin{gather*}
	D(\delta,t_0)=\{(x,y,t)\in \mathcal B_{r}\times\{t\geq t_0\}: \ \ |d_\vartheta({\bf z},t;u_\ast)|\leq \delta\}
\end{gather*}
with some $0<\delta<\varepsilon$ and $t_0\geq t_\ast$. Consider the auxiliary functions 
\begin{gather*}
	U_\vartheta(x,y,t)= t^\vartheta U_1(x,y,t), \quad u(x,y,t)\equiv U_\vartheta(x,y,t)-u_\ast.
\end{gather*}
From \eqref{VN} it follows that there exists $M_1>0$ such that 
\begin{gather}\label{axest}
|d_\vartheta({\bf z},t;u_\ast)|-M_1 t^{-\frac{1}{q}} \leq |u(x,y,t)|\leq |d_\vartheta({\bf z},t;u_\ast)|+M_1 t^{-\frac{1}{q}}
\end{gather}
for all $(x,y,t)\in D(\delta,t_0)$. 
It can easily be checked that
\begin{gather}\label{LUvartheta}
\begin{split}
LU_\vartheta(x,y,t)\equiv & \vartheta t^{-1} U_\vartheta(x,y,t)+t^\vartheta L U_1(x,y,t)\\
	 = &t^{-\frac{n+l}{q}}  U_\vartheta(x,y,t)\left( \lambda_{n+l}+\delta_{n+l,q}\vartheta+\lambda_{n,m} \big(U_\vartheta(x,y,t)\big)^{m-1}+\tilde R(x,y,t) \right)
	\end{split}
\end{gather}
for all $t\geq t_\ast$ and $(x,y)\in \mathcal D(E_0)$ with $\tilde R(x,y,t)=\mathcal O(t^{-\varkappa/q})$, $\varkappa=\min\{1,q\vartheta\}$ as $t\to\infty$. Note that there exists $M_2>0$ such that 
\begin{gather*}
	U_\vartheta(x,y,t)|\tilde R(x,y,t)|\leq M_2 t^{-\frac{\varkappa}{q}}\quad \forall\, (x,y,t)\in D(\delta,t_0).
\end{gather*}
Combining this with \eqref{LUvartheta}, we obtain
\begin{gather}\label{Luest}
L|u(x,y,t)|\leq - t^{-1}U_\vartheta(x,y,t)  |\lambda_{n,m}|\big| (U_\vartheta(x,y,t))^{m-1}-u_\ast^{m-1}\big|+ M_2 t^{-1-\frac{\varkappa}{q}}
\end{gather}
for all $(x,y,t)\in D(\delta,t_0)$.

Consider
\begin{gather*}
	U_3(x,y,t)=|u(x,y,t)| +\left(M_1+ \frac{q}{\varkappa} M_2\right) t^{-\frac{\varkappa}{q}}
	\end{gather*}
 as a Lyapunov function candidate for system \eqref{FulSys}. Taking into account \eqref{axest} and \eqref{Luest}, we get
\begin{gather}\label{LU3est}
	\begin{split}
	&|d_\vartheta({\bf z},t;u_\ast)|\leq U_3(x,y,t)\leq |d_\vartheta({\bf z},t;u_\ast)|+ 2 t^{-\frac{\varkappa}{q}} \left(M_1+ \frac{q}{\varkappa} M_2\right),\\
	& LU_3(x,y,t)\leq - t^{-1}U_\vartheta(x,y,t) |\lambda_{n,m}|\big| (U_\vartheta(x,y,t))^{m-1}-u_\ast^{m-1}\big|-\frac{\varkappa}{q} M_1t^{-1-\frac{\varkappa}{q}}\leq 0
	\end{split}
\end{gather}
 for all $(x,y,t)\in D(\delta,t_0)$.  

Let ${\bf z}(t)\equiv (x(t),y(t))^T$ be a solution of system \eqref{FulSys} with initial data ${\bf z}_0=(x(t_0),y(t_0))^T$ such that $|d_\vartheta({\bf z}_0,t_0;u_\ast)|\leq \delta$ and $\tau_{D}$ be the first exit time of $({\bf z}(t),t)$ from the domain $ D(\varepsilon,t_0)$. Define the function $\tau_D(t)=\min\{\tau_D,t\}$, then ${\bf z}(\tau_D(t))$ is the process stopped at the first exit time. It follows from \eqref{LU3est} that $U_3(x(\tau_{D}(t)),y(\tau_{D}(t)),\tau_{D}(t))$ is a non-negative supermartingale. In this case, the following estimates hold:
\begin{align*}
\mathbb P\left(\sup_{t\geq t_0} |d_\vartheta({\bf z}(t),t;u_\ast)| \geq \varepsilon\right)
& =	\mathbb P\left(\sup_{t\geq t_0}|d_\vartheta({\bf z}(\tau_{D}(t)),\tau_{D}(t);u_\ast)|\geq \varepsilon\right)\\
& \leq \mathbb P\left(\sup_{t\geq t_0}U_3(x(\tau_{D}(t)),y(\tau_{D}(t)),\tau_{D}(t))\geq \varepsilon\right)\leq \frac{U_3(x(t_0),y(t_0),t_0)}{\varepsilon}.
\end{align*}
From \eqref{LU3est} it follows that $U_3(x(t_0),y(t_0),t_0)\leq \delta+2t_0^{-\varkappa/q}(M_1+qM_2/\varkappa)$. Hence, taking 
\begin{gather*}
\delta=\frac{\varepsilon \eta}{2}, \quad  
t_0=\max\left\{t_\ast,\Big(\frac{\varkappa M_1+qM_2}{\varepsilon \eta \varkappa}\Big)^{\frac{q}{\varkappa}}\right\},
\end{gather*} 
we obtain \eqref{defineq2}. 
\qed

\subsection*{Proof of Theorem~\ref{Th5}}
The proof is similar to that of Theorem~\ref{Th3} using the Lyapunov function \eqref{UnFun} and the auxiliary function \eqref{AuxFun} with $N=n+l$. \qed

\section{Stability analysis in case \eqref{ass30}}
\label{Sec43}

\subsection*{Proof of Theorem~\ref{Th6}}
Fix the parameters $\varepsilon>0$ and $\eta>0$, and define
\begin{gather*}
	D(\delta_0,t_\ast)=\{(x,y,t)\in \mathcal B_{r}\times\{t\geq t_\ast\}: \ \ |d_0({\bf z},t;c)|\leq \delta_0\}
\end{gather*}
with a sufficiently small $\delta_0>0$ such that $c+\delta_0<E_0$. Consider the auxiliary functions 
\begin{gather*}
U(x,y,t)= V_{N}(I(x,y),\Phi(x,y),t), \quad u(x,y,t)=  U(x,y,t)-c
\end{gather*}
with $N=n$.
It can easily be checked that
\begin{gather}\label{Lu}
\begin{split}
L|u(x,y,t)|\equiv &{\hbox{\rm sgn}}(u(x,y,t))\mathcal L V_n(E,\varphi,t) \big|_{E=I(x,y), \varphi=\Phi(x,y)}\\
\equiv  & t^{-1} {\hbox{\rm sgn}}(u(x,y,t))\Lambda_n(U_n(x,y,t))+\tilde R_n(x,y,t), 
\end{split}
\end{gather}
for all $t\geq t_\ast$ and $(x,y)\in\mathcal D(E_0)$, where $\tilde R_n(x,y,t)\equiv {\hbox{\rm sgn}}(u(x,y,t))\tilde \Lambda_n(I(x,y),\Phi(x,y),t)$. Note that there exist $0<\delta_1\leq \delta_0$ and $t_1\geq t_\ast$ such that
\begin{gather*}
{\hbox{\rm sgn}}(u(x,y,t))\Lambda_n(U_n(x,y,t))\equiv {\hbox{\rm sgn}}(u(x,y,t))\Big(\Lambda_n(U_n(x,y,t))-\Lambda_n(c)\Big)
\leq -\frac{|\Lambda_n'(c)|}{2} |u(x,y,t)|
\end{gather*}
for all $(x,y,t)\in D(\delta_1,t_1)$. Moreover, from \eqref{RemEst} and \eqref{VN} it follows that 
\begin{gather}\label{axest2}
|d_0({\bf z},t;c)|-M_1 t^{-\frac{1}{q}} \leq |u(x,y,t)|\leq |d_0({\bf z},t;c)|+M_1 t^{-\frac{1}{q}}, \quad 
|\tilde R_n(x,y,t)|\leq M_2 t^{-1-\frac{1}{q}}
\end{gather}
for all $(x,y,t)\in D(\delta_1,t_1)$ with some constants $M_1>0$ and $M_2>0$. 

Consider
\begin{gather*}
	U_0(x,y,t)=|u(x,y,t)|+(M_1+ q M_2) t^{-\frac{1}{q}}
	\end{gather*}
 as a Lyapunov function candidate for system \eqref{FulSys}. From \eqref{Lu} and \eqref{axest2} it follows that  
\begin{gather}\label{LUest2}
	\begin{split}
	&|d_0({\bf z},t;c)|\leq U_0(x,y,t)\leq |d_0({\bf z},t;c)|+ 2 t^{-\frac{1}{q}} (M_1+ q M_2),\\
	& LU_0(x,y,t)\leq -t^{-1}\frac{|\Lambda_n'(c)|}{2} |u(x,y,t)|-\frac{M_1}{q}t^{-1-\frac{1}{q}}\leq 0
	\end{split}
\end{gather}
 for all $(x,y,t)\in D(\delta_1,t_1)$. 

Let ${\bf z}(t)\equiv (x(t),y(t))^T$ be a solution of system \eqref{FulSys} with initial data ${\bf z}_0=(x(t_0),y(t_0))^T$ such that $|d_0({\bf z}_0,t_0;c)|\leq \delta$ with some $0<\delta\leq \delta_1$ and $t_0\geq t_1$. Let $\tau_{D}$ be the first exit time of $({\bf z}(t),t)$ from the domain $D(\varepsilon,t_0)$.
Define the function $\tau_{D}(t)=\min\{\tau_{D},t\}$, then ${\bf z}(\tau_D(t))$ is the process stopped at the first exit time. It follows from \eqref{LUest2} that $U_0(x(\tau_{D}(t)),y(\tau_{D}(t)),\tau_{D}(t))$ is a non-negative supermartingale. Hence, the following estimates hold:
\begin{align*}
\mathbb P\left(\sup_{t\geq t_0} |d_0({\bf z}(t),t;c)| \geq \varepsilon\right)
& =	\mathbb P\left(\sup_{t\geq t_0}|d_0({\bf z}(\tau_{D}(t)),\tau_{D}(t);c)|\geq \varepsilon\right)\\
& \leq \mathbb P\left(\sup_{t\geq t_0}U_0(x(\tau_{D}(t)),y(\tau_{D}(t)),\tau_{D}(t))\geq \varepsilon\right) \leq \frac{U_0(x(t_0),y(t_0),t_0)}{\varepsilon}.
\end{align*}
Since $U_0(x(t_0),y(t_0),t_0)\leq \delta_0+2t_0^{-1/q}(M_1+qM_2)$, we can take 
\begin{gather*}
\delta_0=\min\left\{\delta_1,\frac{\varepsilon \eta}{2}\right\}, \quad  
t_0=\max\left\{t_1,\Big(\frac{M_1+qM_2}{\varepsilon \eta}\Big)^q\right\}
\end{gather*} 
to obtain \eqref{defineq3}. 
\qed

\section{Examples}\label{SecEx}

\subsection*{Example 1} First, consider decaying perturbations of system \eqref{LimSys} with 
\begin{gather}\label{Hexample}
	H_0(x, y)\equiv 1-\cos x + \frac{y^2}{2}.
\end{gather}
It is readily seen that, in this case, the equilibrium $(0, 0)$ of the limiting system \eqref{LimSys} is a center, and the level lines $H_0(x,y)\equiv E$ with $E \in (0, 2)$, lying in the neighborhood of the equilibrium, correspond to $T(E)$-periodic solutions with $\nu(E)\equiv2\pi/T(E)=1-E/8+\mathcal O(E^2)$ as $E\to 0$.

Consider the perturbed system in the following form:
\begin{gather}\label{Ex1}
\begin{split}
&dx=\partial_y H_0(x,y)\,dt, \\ 
&dy=\left(-\partial_x H_0(x,y)+t^{-\frac{h}{q}}\lambda y\right)\,dt+t^{-\frac{p}{q}} \mu \sin x \,dw_2(t),
\end{split}
\end{gather}
with $0<h,p\leq q$ and $\lambda,\mu={\hbox{\rm const}}$. Note that this system is of the form \eqref{FulSys} with $F(x,y,t)\equiv t^{-h/q}\lambda y$, $B_{1,1}(x,y,t)\equiv B_{1,2}(x,y,t)\equiv B_{2,1}(x,y,t)\equiv 0$ and $B_{2,2}(x,y,t)\equiv t^{-p/q}\mu \sin x$. 

Let $h=p=1$ and $q=2$. Then, the changes of the variables described in Section~\ref{Sec3} with $N=1$ and
\begin{gather*}
v_1(E,\varphi)= -\frac{\lambda }{\nu(E)}\int\limits_0^\varphi \left\{\big(Y(\sigma,E)\big)^2  \right\}_\sigma\,d\sigma+\frac{\lambda}{\nu(E)}\left\langle\int\limits_0^\varphi \left\{\big(Y(\sigma,E)\big)^2\right\}_\sigma\,d\sigma\right\rangle_\varphi
\end{gather*}
transform \eqref{Ex1} to \eqref{Veq} with
\begin{gather*}
\Lambda_1(v)\equiv \frac{\lambda}{\nu(v)} \left\langle\big(Y(\varphi,v)\big)^2\right\rangle_\varphi=v\big(\lambda+\mathcal O(v)\big), \quad v\to 0.
\end{gather*}
It can easily be checked that the transformed system satisfies \eqref{ass1}, \eqref{ass21} and \eqref{ass3} with $n=1$, $\lambda_n=\lambda$, and $\sigma=1>n/q$. Hence, it follows from Theorem~\ref{Th2} that the equilibrium $(0,0)$ of system \eqref{Ex1} is exponentially stable in probability if $\lambda<0$. From Theorem~\ref{Th3} it follows that  the equilibrium $(0,0)$ is unstable in probability if $\lambda>0$. In this case the stability of the equilibrium is determined by the sign of the coefficient $\lambda$ as in the deterministic system with $\mu=0$ (see Fig.~\ref{FigEx11}).

\begin{figure}
\centering
\subfigure[$\mu=0$ ]{\includegraphics[width=0.4\linewidth]{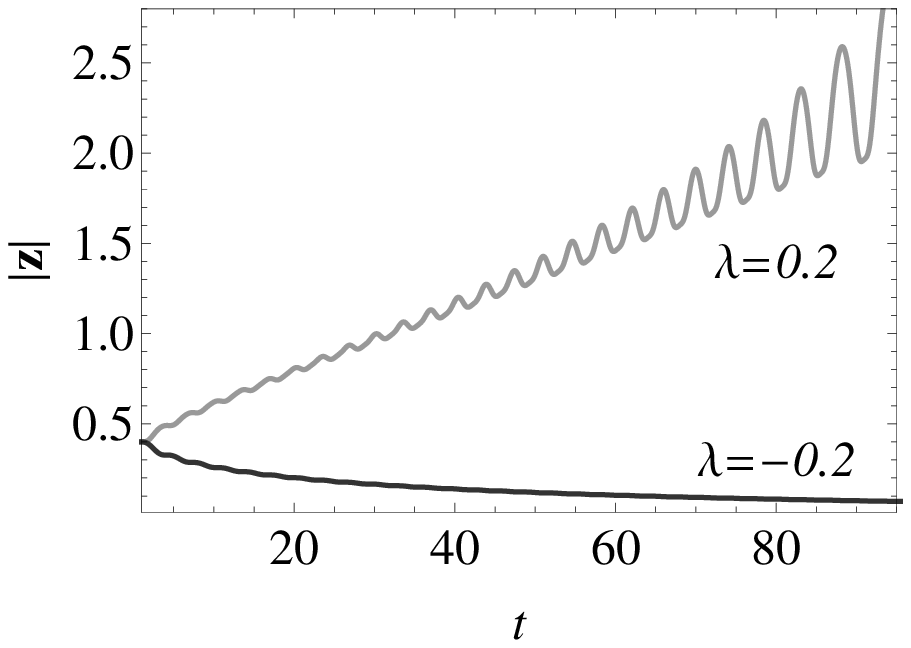}}
\hspace{4ex}
\subfigure[$\mu=1$]{\includegraphics[width=0.4\linewidth]{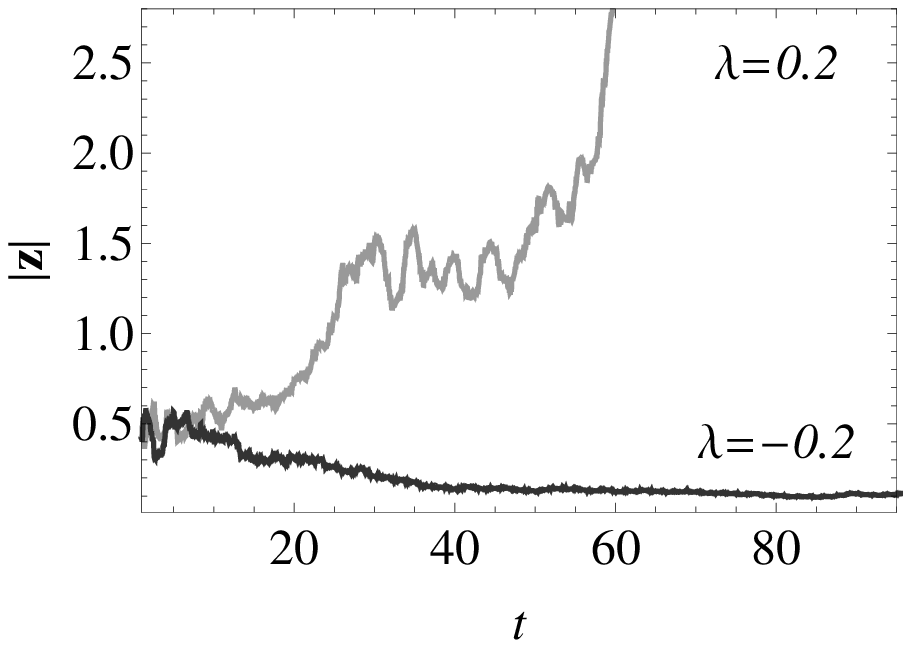}}
\caption{\small The evolution of $|{\bf z}(t)|=\sqrt{x^2(t)+y^2(t)}$ for sample paths of the solutions to system \eqref{Ex1} with $h=p=1$, $q=2$ and initial data $x(1)=0.4$, $y(1)=0$.} \label{FigEx11}
\end{figure}

Let $h=q=2$ and $p=1$. Under the transformation described in Section~\ref{Sec3} with $N=2$, $v_1(E,\varphi)\equiv 0$ and
\begin{align*}
v_2(E,\varphi)\equiv & -\frac{1}{\nu(E)}\int\limits_0^\varphi \left\{\lambda \big(Y(\sigma,E)\big)^2 +\frac{\mu^2}{2} \big(\sin (X(\sigma,E))\big)^2 \right\}_\sigma\,d\sigma \\
& +\frac{1}{\nu(E)}\left\langle\int\limits_0^\varphi \left\{\lambda \big(Y(\sigma,E)\big)^2 +\frac{\mu^2}{2} \big(\sin (X(\sigma,E))\big)^2\right\}_\sigma\,d\sigma\right\rangle_\varphi
\end{align*}
system \eqref{Ex1} is reduced to \eqref{Veq} with $\Lambda_1(v)\equiv 0$ and 
\begin{align*}
\Lambda_2(v)\equiv& \frac{1}{\nu(v)} \left\langle\lambda \big(Y(\sigma,v)\big)^2 +\frac{\mu^2}{2} \big(\sin (X(\sigma,v))\big)^2  \right\rangle_\varphi =v\left(\lambda+\frac{\mu^2}{2}+\mathcal O(v)\right), \quad v\to 0.
\end{align*}
It is readily seen that the transformed system satisfies \eqref{ass1}, \eqref{ass21} and \eqref{ass3} with $n=2$, $\lambda_n=\lambda+\mu^2/2$, and $\sigma=1=n/q$. Applying Theorem~\ref{Th2}, we see that if $\lambda<-\mu^2/2$, then for all $\kappa\in (0,1)$ the equilibrium $(0,0)$ is stable in probability with the weight $t^{(1-\kappa)|\lambda_n|/2}$. From Theorem~\ref{Th3} it follows that the equilibrium is unstable if $\lambda>0$, and unstable with a weight if $-\mu^2/2< \lambda\leq 0$ (see Fig.~\ref{FigEx12}). In the latter case, from Theorem~\ref{AsL} it follows that the equilibrium $(0,0)$ is practically stable.

\begin{figure} 
\centering
\subfigure[$\mu=0$ ]{\includegraphics[width=0.4\linewidth]{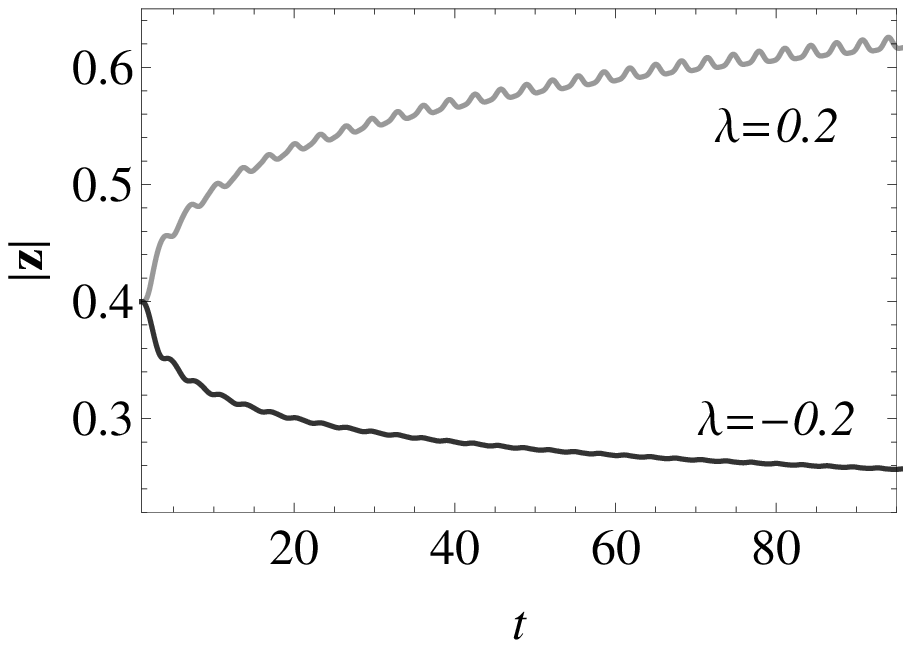}}
\hspace{4ex}
\subfigure[$\mu=1$]{\includegraphics[width=0.4\linewidth]{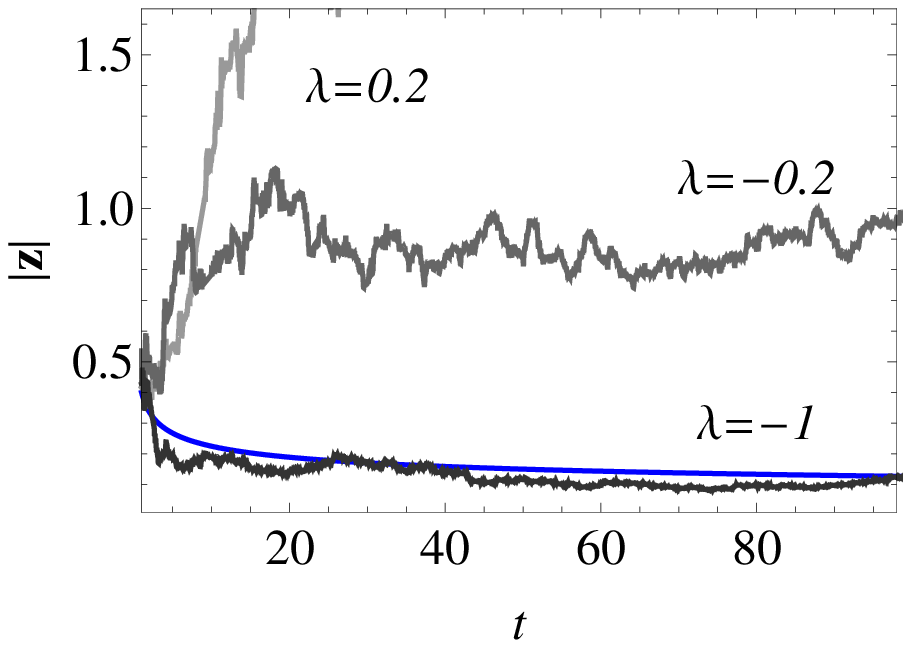}}
\caption{\small The evolution of $|{\bf z}(t)|=\sqrt{x^2(t)+y^2(t)}$ for sample paths of the solutions to system \eqref{Ex1} with $h=q=2$, $p=1$ and initial data $x(1)=0.4$, $y(1)=0$. The blue curve corresponds to $|{\bf z}|=t^{-1/4}|{\bf z}(1)|$.} \label{FigEx12}
\end{figure}

\subsection*{Example 2} Consider the perturbed system
\begin{gather}\label{Ex2}
\begin{split}
&dx=\partial_y H_0(x,y)\,dt, \\ 
&dy=\left(-\partial_x H_0(x,y)+t^{-\frac{1}{2}}F_2(x,y)+t^{-1}F_4(y)\right)\,dt+\left(t^{-\frac{1}{4}} B_{2,2,1}(x,y)+t^{-\frac{1}{2}}B_{2,2,2}(x) \right)\,dw_2(t),
\end{split}
\end{gather}
with $H_0(x,y)$ defined by \eqref{Hexample}, 
\begin{gather*}
F_2(x,y)= \frac{a_2 x^2y}{1+x^2}, \quad F_4(y)=a_4 y, \quad B_{2,2,1}(x,y)=\frac{b_1 xy}{\sqrt{1+x^2}}, \quad B_{2,2,2}(x)=b_2 x,
\end{gather*}
where $a_2,a_4,b_1,b_2={\hbox{\rm const}}$. It is readily seen that this system is of the form \eqref{FulSys} with $q=4$.
The changes of the variables described in Section~\ref{Sec3} with $N=4$, $v_1(E,\varphi)\equiv 0$, 
\begin{align*}
v_2(E,\varphi)\equiv & -\frac{1}{2\nu(E)}\int\limits_0^\varphi \left\{ 2YF_2(X,Y) + \big(B_{2,2,1}(X,Y)\big)^2 \right\}_\sigma\,d\sigma \\
& +\frac{1}{2\nu(E)}\left\langle\int\limits_0^\varphi \left\{ 2YF_2(X,Y) + \big(B_{2,2,1}(X,Y)\big)^2 \right\}_\sigma\,d\sigma\right\rangle_\varphi,\\
v_3(E,\varphi)\equiv & -\frac{1}{\nu(E)}\int\limits_0^\varphi \left\{  B_{2,2,1}(X,Y)B_{2,2,2}(X)  \right\}_\sigma\,d\sigma  +\frac{1}{\nu(E)}\left\langle\int\limits_0^\varphi \left\{  B_{2,2,1}(X,Y)B_{2,2,2}(X)\right\}_\sigma\,d\sigma\right\rangle_\varphi,\\
v_4(E,\varphi)\equiv & -\frac{1}{2\nu(E)}\int\limits_0^\varphi \left\{ 2Y F_4(Y)+  \big(B_{2,2,2}(X)\big)^2 +2R_4(E,\sigma) \right\}_\sigma\,d\sigma \\
& +\frac{1}{2\nu(E)}\left\langle\int\limits_0^\varphi \left\{ 2 Y F_4(Y)+  \big(B_{2,2,2}(X)\big)^2 +2R_4(E,\sigma)\right\}_\sigma\,d\sigma\right\rangle_\varphi
\end{align*}
reduce \eqref{Ex2} to \eqref{Veq} with $\Lambda_1(E)\equiv 0$,
\begin{align*}
&\Lambda_2(v)\equiv \frac{1}{2\nu(v)} \left\langle 2YF_2(X,Y) + \big(B_{2,2,1}(X,Y)\big)^2\right\rangle_\varphi=\frac{v^2}{4}\left(2a_2+b_1^2+\mathcal O(v)\right),\\
&\Lambda_3(v)\equiv  \frac{1}{\nu(v)} \left\langle  B_{2,2,1}(X,Y)B_{2,2,2}(X) \right\rangle_\varphi=\mathcal O(v^3),\\
&\Lambda_4(v)\equiv  \frac{1}{2\nu(v)} \left\langle 2 Y F_4(Y)+  \big(B_{2,2,2}(X)\big)^2 +2R_4(v,\varphi) \right\rangle_\varphi= \frac{v}{2}\left(2a_4+b_2^2+\mathcal O(v)\right)
\end{align*}
as $v\to 0$,
where 
\begin{align*}
R_4(E,\varphi) \equiv &  -\partial_E\Lambda_2(E)v_2(E,\varphi)   +F_2(X ,Y ) \Big( Y\partial_E+\partial_y \Phi(X ,Y )\partial_\varphi \Big) v_2(E,\varphi)\\
	& +\frac{1}{2}\big(B_{2,2,1}(X ,Y )\big)^2\Big( \partial_E +  \partial_y^2 \Phi(X ,Y )\partial_\varphi\Big)v_2(E,\varphi)=\mathcal O(E^3)
\end{align*}
as $E\to 0$ uniformly for all $\varphi\in \mathbb R$. We see that the transformed system satisfies \eqref{ass1} and \eqref{ass22} with $n=m=l=2$, $\lambda_{n,m}=(2a_2+b_1^2)/4$ and $\lambda_{n+l}=(2a_4+b_2^2)/2$. In this case, $n+l=q$, $\gamma_{n+l}(t)\equiv t$, $\vartheta=1/2$, $u_\ast=2|2a_4+b_2^2+1|/|2a_2+b_1^2|$ and $d_\vartheta({\bf z},t;u_\ast)\equiv t^{1/2}H_0(x,y)-u_\ast$. It follows from Theorem~\ref{Th4} that the equilibrium $(0,0)$ of system \eqref{Ex2} is stable in probability if $a_2<-b_1^2/2$ and $a_4<-b_2^2/2$. Moreover, for all $\kappa\in(0,1)$ the equilibrium $(0,0)$ is stable with the weight $t^{(\vartheta +(1-\kappa)|\lambda_{n+l}+\vartheta|)/2}$ if $a_4<-(b_2^2+1)/2$ (see Fig.~\ref{FigEx21}, a).  If $a_2<-b_1^2/2$ and $a_4>-b_2^2/2$, the equilibrium $(0,0)$ is stable in the sense of \eqref{defineq2}. In the latter case, $H_0(x(t),y(t))\sim u_\ast t^{-1/2}$ as $t\to\infty$ with high probability for solutions starting in the vicinity of the equilibrium (see Fig.~\ref{FigEx22}).

\begin{figure} 
\centering
\subfigure[$a_2=1$, $a_4=-5/4$, $b_1=4$, $b_2=1$ ]{\includegraphics[width=0.4\linewidth]{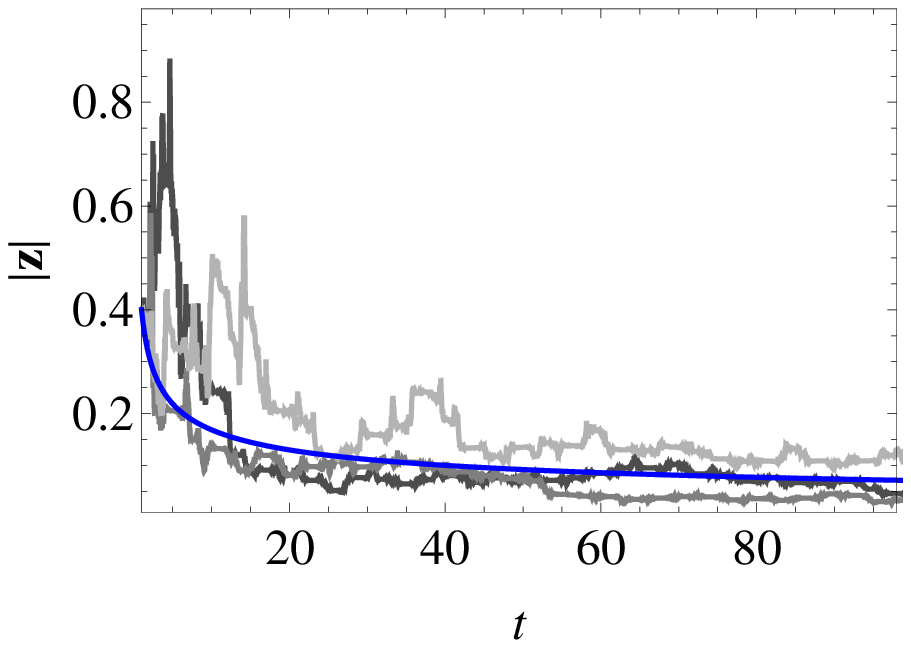}}
\hspace{4ex}
\subfigure[$a_2=a_4=0.1$, $b_1=0$, $b_2=1$]{\includegraphics[width=0.4\linewidth]{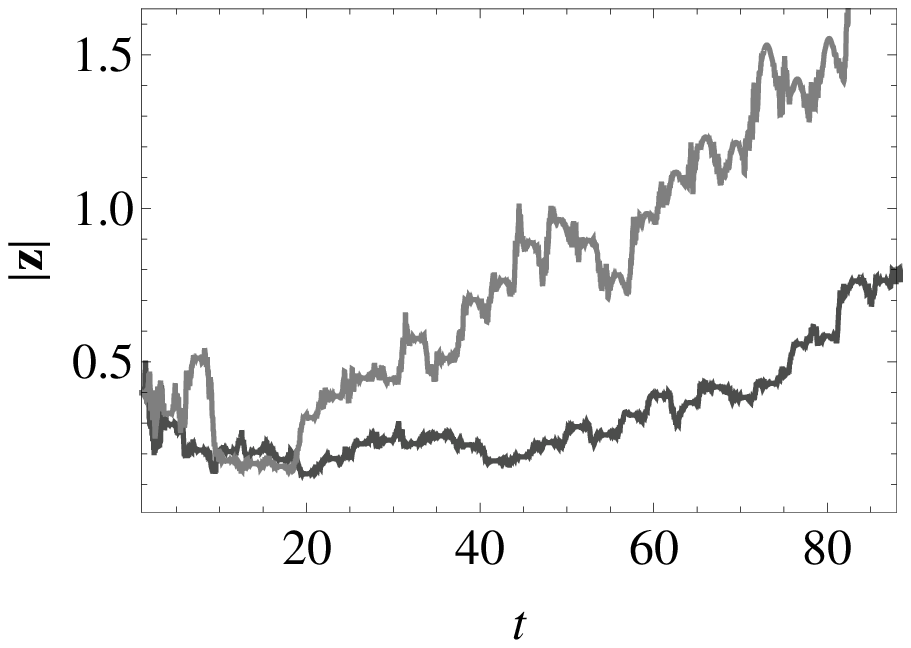}}
\caption{\small The evolution of $|{\bf z}(t)|=\sqrt{x^2(t)+y^2(t)}$ for sample paths of the solution to system \eqref{Ex2} with initial data $x(1)=0.4$, $y(1)=0$ and different values of the parameters. The blue curve corresponds to $|{\bf z}(1)| t^{-3/8}$.} \label{FigEx21}
\end{figure}

\begin{figure} 
\centering
{\includegraphics[width=0.4\linewidth]{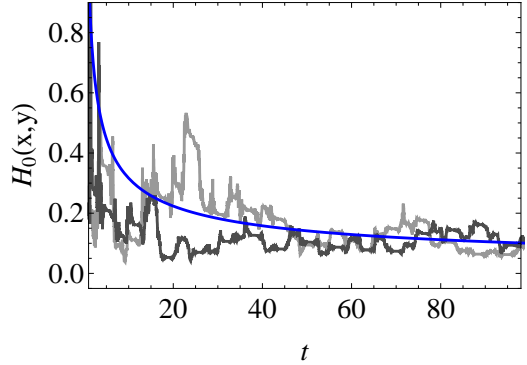}}
\caption{\small The evolution of $H_0(x(t),y(t))$ for sample paths of the solution to system \eqref{Ex2} with $a_2=-2$, $a_4=-1/4$, $b_1=1$, $b_2=1$ and various initial data. The blue curve corresponds to $u_\ast t^{-1/2}$.} \label{FigEx22}
\end{figure}

If $b_1=0$, then system \eqref{Ex2} additionally satisfies \eqref{ass3} with $\sigma=1=(n+l)/q$ and $\mu=b_2^2$. 
From Theorem~\ref{Th5} it follows that the equilibrium $(0,0)$ is unstable in probability if $a_2>0$ and $a_4>0$ (see Fig.~\ref{FigEx21}, b). 

\subsection*{Example 3} Finally, consider the perturbed system
\begin{gather}\label{Ex3}
\begin{split}
 dx=y\,dt, \quad  
 dy=\left(-x+t^{-1}\frac{(a_1+a_2x^2) y}{1+x^2+y^2}\right)\,dt+t^{-\frac{1}{2}} \mu x \,dw_2(t)
\end{split}
\end{gather}
with $a_1,a_2,\mu={\hbox{\rm const}}$. It is clear that this system is of the form \eqref{FulSys} with $q=2$, $H_0(x,y)\equiv|{\bf z}|^2/2$, $\nu(E)\equiv 1$. In this case, $X(\varphi,E)=\sqrt{2E}\cos\varphi$ and $Y(\varphi,E)=-\sqrt{2E}\sin \varphi$. The changes of the variables described in Section~\ref{Sec3} with $N=2$, $v_1(E,\varphi)\equiv 0$ and 
\begin{align*}
v_2(E,\varphi)\equiv & - \int\limits_0^\varphi \left\{  \frac{  (a_1+a_2 X^2)  Y^2}{1+2E}+ \frac{\mu^2 }{2}X^2   \right\}_\sigma\,d\sigma   +\left\langle\int\limits_0^\varphi \left\{  \frac{(a_1+a_2 X^2)  Y^2}{1+2E}+ \frac{\mu^2 }{2} X^2    \right\}_\sigma\,d\sigma\right\rangle_\varphi 
\end{align*}
reduce \eqref{Ex3} to \eqref{Veq} with $\Lambda_1(v)=0$ and
\begin{gather*}
\Lambda_2(v)\equiv \left\langle \frac{  (a_1+a_2 X^2(\varphi,v)) Y^2(\varphi,v)     }{1+2v} + \frac{\mu^2}{2}   X^2(\varphi,v)  \right\rangle_\varphi=\frac{v \big(2a_1+\mu^2+v(a_2+2\mu^2)\big)}{2(1+2v)}.
\end{gather*}
It is readily seen that the transformed system satisfies \eqref{ass1}, \eqref{ass21} and \eqref{ass3} with $n=2$, $\lambda_n=a_1+\mu^2/2$ and $\sigma=n/q=1$. It follows from Theorem~\ref{Th3} that the equilibrium $(0,0)$ is unstable if $a_1>0$. If, in addition, $a_2<-2\mu^2$, then system \eqref{Ex3} satisfies \eqref{ass30} with $c=(2a_1+\mu^2)/|a_2+2\mu^2|>0$ such that $\Lambda_2(c)=0$ and $\Lambda_2'(c)=-(2a_1+\mu^2)/(2(1+2c))<0$. 
Hence, by applying Theorem~\ref{Th6}, we conclude that there exists a stable cycle $|{\bf z}_\ast(t)|\approx \sqrt{2c}$ (see Fig.~\ref{FigEx3}).

\begin{figure}
\centering
{\includegraphics[width=0.4\linewidth]{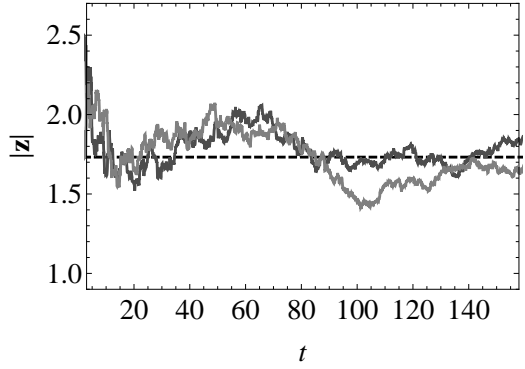}}
\caption{\small The evolution of $|{\bf z}(t)|=\sqrt{x^2(t)+y^2(t)}$ for sample paths of the solutions to system \eqref{Ex3} with $\mu=1/2$, $a_1=1$ and $a_2=-2$. The dashed line corresponds to $|{\bf z}|=	\sqrt{2c}=\sqrt{3}$. } \label{FigEx3}
\end{figure}

\section{Conclusion}

Thus, possible bifurcations associated with changes in stochastic stability of the equilibrium in asymptotically Hamiltonian systems with multiplicative noise have been described. Through a careful nonlinear analysis based on a combination of the averaging method and the construction of stochastic Lyapunov functions we have shown that depending on the structure of decaying perturbations the equilibrium of the limiting system becomes stable (exponentially, polynomially or practically), or loses stability. In the latter case, stable cycles may appear.

Note that for the considered class of systems the linearization method for stability analysis fails. In particular, Example 2 demonstrates this. Indeed, the linearization of system \eqref{Ex2} with $b_1=0$ in the vicinity of the equilibrium $(0,0)$ is given by system \eqref{ex0} with $\lambda=a_4$ and $\mu=b_2$. Applying Theorems~\ref{Th1}, \ref{Th2} and \ref{Th3} shows that if $a_4<-b_2^2/2$, the equilibrium is stable in probability, and if $a_4>0$, the equilibrium is unstable (compare Example 1). However, the equilibrium is stable in the full system if $a_4>-b_2^2/2$ and $a_2<0$ (see Example 2). This contradicts the instability obtained for the corresponding linearization.

Note also that if perturbations do not preserve the equilibrium, the proposed method cannot be applied directly. Bifurcation phenomena in this case should be considered separately. This will be discussed elsewhere.

\section*{Acknowledgments}
Research is supported by the Russian Science Foundation grant 19-71-30002.

}

\end{document}